\documentclass{tran-l}

\newtheorem{theorem}{Theorem}[section]
\newtheorem{lemma}[theorem]{Lemma}
\newtheorem{Cor}{Corollary}

\theoremstyle{definition}

\newtheorem{example}[theorem]{Example}

\theoremstyle{remark}
\newtheorem{remark}[theorem]{Remark}
\theoremstyle{remark}
\newtheorem{remarks}[theorem]{Remarks}

\numberwithin{equation}{section}

\newcommand{\from}[2]{\left|{\displaystyle{\phantom\int}\!\!\!\!\!}_{#1}^{#2}
  \right.}
\newcommand{\Li}{\mathrm{Li}}
\newcommand{\eu}{\stackrel{?}{=}}
\renewcommand{\a}{\alpha}
\renewcommand{\b}{\beta}
\newcommand{\g}{\gamma}
\newcommand{\G}{\Gamma}
\renewcommand{\d}{\delta}
\renewcommand{\l}{\lambda}
\newcommand{\Z}{{\mathbf Z}}
\newcommand{\R}{{\mathbf R}}
\newcommand{\Q}{{\mathbf Q}}
\newcommand{\C}{{\mathbf C}}
\newcommand{\w}{\omega}
\newcommand{\W}{\Omega}
\newcommand{\z}{\zeta}
\newcommand{\iu}{\int_0^1}
\newcommand{\us}{\{1\}}
\renewcommand{\(}{\left(}
\renewcommand{\)}{\right)}

\newcommand{\Cat}[2]{\mathop{\mathrm{\bf Cat}}\limits_{#1}^{#2}}
\def\eps{\varepsilon}

\begin{document}

\title[Special Values of Multiple Polylogarithms]{Special Values of Multiple
Polylogarithms}

% at present the "communicated by" line appears only in ERA and PROC
\commby{Alice Silverberg}

\author[Borwein]{Jonathan M. Borwein}
\address{Centre for Experimental and Constructive Mathematics\\
         Simon Fraser University\\
         Burnaby, B.C., V5A$\;$1S6\\
         Canada}
\email{ jborwein@cecm.sfu.ca}
\thanks{The research of the first author was supported by
        NSERC and the Shrum Endowment of Simon Fraser University.}

\author[Bradley]{David M. Bradley}
\address{University of Maine\\
         Department of Mathematics and Statistics\\
         5752 Neville Hall\\
         Orono, Maine 04469--5752\\
         U.S.A.}
\email{bradley@gauss.umemat.maine.edu, dbradley@member.ams.org}

\author[Broadhurst]{David J. Broadhurst}
\address{Physics Department, Open University\\
         Milton Keynes, MK7$\;$6AA\\
         United Kingdom}
\email{D.Broadhurst@open.ac.uk}

\author[Lison\v ek]{Petr Lison\v ek}
\address{Centre for Experimental and Constructive Mathematics\\
         Simon Fraser University\\
         Burnaby, B.C., V5A$\;$1S6\\
         Canada}
\email{lisonek@cecm.sfu.ca}

%\begin{ao}
%David M. Bradley, University of Maine, Department of Mathematics
%and Statistics, 5752 Neville Hall, Orono, Maine 04469-5752, U.S.A.
%(e-mail: dbradley@member.ams.org, bradley@gauss.umemat.maine.edu,
%dbradley@cecm.sfu.ca).
%\end{ao}

\date{%Submitted
July 10, 1998. Revised August 9, 1999.}

\subjclass{Primary: 40B05, 33E20; Secondary: 11M99, 11Y99}

\keywords{Euler sums, Zagier sums, multiple zeta values,
          polylogarithms, multiple harmonic series,
          quantum field theory, knot theory, Riemann zeta function.}

\begin{abstract}
Historically, the polylogarithm has attracted specialists and
non-specialists alike with its lovely evaluations. Much the same
can be said for Euler sums (or multiple harmonic sums), which,
within the past decade, have arisen in combinatorics, knot theory
and high-energy physics. More recently, we have been forced to
consider multidimensional extensions encompassing the classical
polylogarithm, Euler sums, and the Riemann zeta function.  Here,
we provide a~general framework within which previously isolated
results can now be properly understood. Applying the theory
developed herein, we prove several previously conjectured
evaluations, including an intriguing conjecture of Don Zagier.
\end{abstract}

\maketitle

%\abbreviations{MZV -- multiple zeta value}

%------------------------------------------------------------------------------
% End of journal.top
%------------------------------------------------------------------------------

\section{Introduction}
\label{S1}
We are going to study a~class of multiply nested sums of the form
\begin{equation}
    \l{s_1,\dots,s_k\choose b_1,\dots,b_k}
    :=\sum_{\nu_1,\ldots,\nu_k=1}^\infty \;
    \prod_{j=1}^k b_j^{-\nu_j}
       \bigg(\sum_{i=j}^k\nu_i\bigg)^{-s_j},
    \label{polydef}
\end{equation}
and which we shall refer to as {\em multiple polylogarithms.}
When $k=0$, we define $\l(\{\}):=1$, where $\{\}$ denotes the
empty string. When $k=1$, note that
\begin{equation}
   \l{s\choose b}
    = \sum_{\nu=1}^\infty \frac{1}{{\nu}^s b^{\nu}}
    = {\Li}_s\(\frac1b\)\label{lidef}
\end{equation}
is the usual polylogarithm~\cite{Lewin1,Lewin2}
when $s$ is a~positive integer and $|b|\ge1$.
Of course, the polylogarithm~(\ref{lidef})
reduces to the Riemann zeta function~\cite{Ed,Ivic,Tit}
\begin{equation}
   \z(s) = \sum_{\nu=1}^\infty \frac1{{\nu}^s},\quad \Re(s)>1,
   \label{Riemann}
\end{equation}
when $b=1$.  More generally, for any $k>0$ the substitution
$n_j=\sum_{i=j}^k \nu_i$ shows that our multiple
polylogarithm~(\ref{polydef}) is related to
Goncharov's~\cite{Gonch3} by the equation
\[
   {\Li}_{s_k,\ldots,s_1}(x_k,\ldots,x_1)
   =\l{s_1,\dots,s_k\choose y_1,\dots,y_k},
   \quad \mbox{where} \quad y_j :=\prod_{i=1}^j x_i^{-1},
\]
and
\begin{equation}
   {\Li}_{s_k,\ldots,s_1}(x_k,\ldots,x_1)
   :=\sum_{n_1>\cdots>n_k>0}\;\prod_{j=1}^k n_j^{-s_j} {x_j}^{n_j}.
\label{Li-nest}
\end{equation}
With each $x_j=1$, these latter sums (sometimes called {\em ``Euler sums''}),
have been studied previously at various levels of
generality~\cite{Bail,BBG,BBB,BG,DJB1,BGK,BK2,BK1,Flaj,Hoff1,Hoff2,Moen,Mar,Niels},
the case $k=2$ going back to Euler~\cite{LE}.
Recently, Euler sums have arisen in combinatorics
(analysis of quad-trees \cite{FLLS,LL}
and of lattice reduction algorithms \cite{DFV}),
knot theory~\cite{BGK,BK2,BK1,LeM},
and high-energy particle physics~\cite{DJB1} (quantum field theory).
There is also quite sophisticated work relating polylogarithms and their
generalizations to arithmetic and algebraic geometry, and
to algebraic
$K$-theory~\cite{Beil,Brow2,Brow,Gonch1,Gonch2,Gonch3,Woj3,Woj2,Woj}.

In view of these recent applications and the well-known fact that
the classical polylogarithm~(\ref{lidef}) often arises in physical problems
via the multiple integration of rational forms, one might expect
that the more general multiple polylogarithm~(\ref{polydef}) would
likewise find application in a~wide variety of physical contexts.
Nevertheless, lest it be suspected that the authors have embarked
on a~program of generalization for its own sake, let the reader be
assured that it was only with the greatest reluctance that we
arrived at the definition~(\ref{polydef}). On the one hand, the
polylogarithm~(\ref{lidef}) has traditionally been studied as
a~function of $b$ with the positive integer $s$ fixed; while on
the other hand, the study of Euler sums has almost exclusively
focused on specializations of the nested sum~(\ref{Li-nest}) in which
each $x_j=\pm1.$   However, we have found,
in the course of our investigations, that a~great deal of insight
is lost by ignoring the interplay between these related sums when
both sequences of parameters are permitted to
vary. Indeed, it is our view that it is {\it impossible} to fully
understand the sums~(\ref{lidef}--\ref{Li-nest}) without viewing
them as members of a~broader class of multiple
polylogarithms.

That said, one might legitimately ask why we chose to adopt the
notation~(\ref{polydef}) in favour of Goncharov's~(\ref{Li-nest}),
inasmuch as the latter is a direct generalization
of the $\Li_n$ notation for the classical polylogarithm.  As a
matter of fact, the notation~(\ref{Li-nest}) (with argument list
reversed) was our original choice.  However, as we reluctantly
discovered, it turns out that the
notation~(\ref{polydef}), in which the second row of parameters
comprises the reciprocated running product of the argument list
in~(\ref{Li-nest}),
is more suitable for our purposes here.  In particular, our
``running product'' notation~(\ref{polydef}), in addition to
simplifying the iterated integral representation~(\ref{iterint})
(cf.~\cite{Gonch1} Theorem 16) and the various duality formulae
(Section~\ref{S6}---see
eg.\ equations~(\ref{deltadualsym}) and (\ref{deltadual})),
brings out much more clearly the relationship (Subsection~\ref{S52})
between the partition integral (Subection~\ref{S41}),
in which running products necessarily arise in the integrand;
and ``stuffles'' (Subsections~\ref{S51}, \ref{S51-half}).
It seems also that boundary cases
of certain formulae for alternating sums must be
treated separately unless running product notation is used.
Theorem~\ref{T5} with $n=0$ (Section~\ref{S8}) provides an
example of this.

Don Zagier (see eg.~\cite{Zag}) has argued persuasively in favour
of studying special values of zeta functions at integer arguments,
as these values ``often seem to dictate the most important
properties of the objects to which the zeta functions are
associated.''  It seems appropriate, therefore, to focus on the
values the multiple polylogarithms~(\ref{polydef}) take when the
$s_j$ are restricted to the set of positive integers, despite the
fact that the sums~(\ref{polydef}) and their special cases have
a~rich structure as analytic functions of the complex variables
$s_j$.  However, we allow the parameters $b_j$ to take on complex
values, with each $|b_j|\ge 1$ and $(b_1,s_1)\ne(1,1)$ to ensure
convergence.

Their importance notwithstanding, we feel obliged to confess that our
interest in special values extends beyond mere utilitarian concerns.
Lewin~\cite{Lewin1} (p.~xi) writes of a~``school-boy fascination''
with certain numerical results, an attitude which we
whole-heartedly share.  In the hope that the reader might also
be convinced of the intrinsic beauty of the subject, we offer two
modest examples.   The first~\cite{Hoff1,LeM},
\[
  \sum_{\nu_1,\ldots,\nu_k=1}^\infty \;
  \prod_{j=1}^k \frac1{(\nu_j+\cdots+\nu_k)^2}
    =\frac{\pi^{2k}}{(2k+1)!},\qquad 0\le k\in\Z,
\]
generalizes Euler's celebrated result
\[
   \z(2) = \sum_{\nu=1}^\infty\frac1{\nu^2}=\frac{\pi^2}6,
\]
and is extended to all even positive integer arguments in~\cite{BBB}.
The second (see Corollary~\ref{C1} of Section~\ref{S8}),
\[
\begin{split}
   \sum_{\nu_1,\ldots,\nu_k=1}^\infty \;
   \prod_{j=1}^k\frac{(-1)^{\nu_j+1}}{\nu_j+\cdots+\nu_k}
   &=\sum_{\nu_1,\ldots,\nu_k=1}^\infty \;
   \prod_{j=1}^k\frac{1}{2^{\nu_j}(\nu_j+\cdots+\nu_k)}\\
   &=\frac{\(\log 2\)^k}{k!},\qquad 0\le k\in\Z
   \label{Mac}
\end{split}
\]
can be viewed as a~multidimensional extension of the
elementary ``dual'' Maclaurin series evaluations
\[
   \sum_{\nu=1}^\infty\frac{(-1)^{\nu+1}}\nu
   =\sum_{\nu=1}^\infty\frac1{\nu 2^{\nu}} =\log 2,
\] and leads to deeper questions of duality (Section~\ref{S6}) and
computational issues related to series acceleration
(Section~\ref{S7}).   We state additional results in the next
section and outline connections to combinatorics and $q$-series.
In Section~\ref{S4}, we develop several different integral
representations, which are then used in subsequent sections to
classify various types of identities that multiple polylogarithms
satisfy.  Sections~\ref{S8} through~\ref{S10} conclude the paper
with proofs of previously conjectured evaluations, including an
intriguing conjecture of Zagier~\cite{Zag} and its generalization.

\section{Definitions and Additional Examples}
\label{S2}

A useful specialization of the general multiple
polylogarithm~(\ref{polydef}), which is at the same time an
extension of the polylogarithm~(\ref{lidef}), is the case in which
each $b_j=b$. Under these circumstances, we write
\begin{equation}
   \l_b(s_1,\dots,s_k):=
   \l{s_1,\dots,s_k\choose b,\dots,b}
   =\sum_{\nu_1,\ldots,\nu_k=1}^\infty \;\prod_{j=1}^k b^{-\nu_j}
      \bigg(\sum_{i=j}^k\nu_i\bigg)^{-s_j},
   \label{polysubdef}
\end{equation}
and distinguish the cases $b=1$ and $b=2$ with special symbols:
\begin{equation}
\label{def-zeta-lambda}
\z:=\l_1 \mbox{\ \ and\ \ } \d:=\l_2.
\end{equation}
The latter $\d$-function represents an iterated sum
extension of the polylogarithm~(\ref{lidef}) with argument one-half,
and will play a~crucial role in computational issues (Section~\ref{S7}) and
``duality'' identities such as~(\ref{Mac}).
The former coincides with~(\ref{Li-nest}) when $k>0$, each $x_j=1$,
and the order of the argument list is reversed,
and hence can be viewed as a~multidimensional unsigned Euler sum.
We will follow Zagier~\cite{Zag} in referring to these as
{\em ``multiple zeta values''} or {\em ``MZVs''} for short.
By specifying each $b_j=\pm 1$ in~(\ref{polydef}), {\em alternating
Euler sums}~\cite{BBB} are recovered, and in this case, it is convenient to
combine the strings of exponents and signs into a~single string with
$s_j$ in the $j$th position when $b_j=+1$, and $s_j-$ in the
$j$th position when $b_j=-1.$
To avoid confusion, it should be also noted that in \cite{BBB}
the alternating Euler sums were studied using the notation
\[
\z(s_1,\ldots,s_k)
   :=\sum_{n_1>\cdots>n_k>0}\;
    \prod_{j=1}^k n_j^{-|s_j|}{\sigma_j^{-n_j}}
\]
where $s_1,\ldots,s_k$ are non-zero integers
and $\sigma_j:=\mbox{signum}(s_j)$.

Additionally, $n$ repetitions
of a~substring $U$ will be denoted by $U^n$.  Thus, for example,
\[
   \l(\{2-,1\}^n) := \l{2,1,\dots,2,1\choose -1,1,\dots,-1,1}
   =\sum_{\nu_1,\ldots,\nu_{2n}=1}^\infty \;\prod_{j=1}^n
   \frac{(-1)^{\nu_{2j-1}}}{\(\sum_{i=2j-1}^k\nu_i\)^2
   \(\sum_{i=2j}^k\nu_i\)}.
\]

{\em Unit Euler sums,} that is those sums~(\ref{polydef}) in which each
$s_j=1$, are also important enough to be given a~distinctive notation.
Accordingly, we define
\begin{equation}
   \mu(b_1,\dots,b_k) := \l{1,\dots,1\choose b_1,\dots,b_k}
   =\sum_{\nu_1,\ldots,\nu_k=1}^\infty \;\prod_{j=1}^k b_j^{-\nu_j}
      \bigg(\sum_{i=j}^k\nu_i\bigg)^{-1}.
   \label{unitEuler}
\end{equation}

To entice the reader, we offer a~small but representative
sample of evaluations below.

%{\sc Example 2.1.}
\begin{example}\label{ex:z21}
Euler showed that
\[
   \z(2,1) = \sum_{n=1}^\infty\frac{1}{n^2}\sum_{k=1}^{n-1}\frac1{k}
   =\sum_{n=1}^\infty\frac1{n^3} = \z(3),
\]
and more generally~\cite{LE,Niels}, that
\[
   2\z(m,1) = m\z(m+1) - \sum_{k=1}^{m-2}\z(m-k)\z(k+1),
   \quad 2\le m\in\Z.
\]
The continued interest in Euler sums is evidenced by the fact
that a~recent American Mathematical Monthly problem \cite{monthly}
effectively asks for the proof of $\z(2,1)=\z(3)$.
\end{example}

Two examples of non-alternating, arbitrary depth evaluations for
all nonnegative integers $n$ are provided by
\begin{example}
%{\sc Example 2.2.}
\[
   \z(\{3,1\}^n) = 4^{-n}\z(\{4\}^n) = \frac{2\pi^{4n}}{(4n+2)!},
\]
previously conjectured by Don Zagier~\cite{Zag} and proved
herein (see Section~\ref{S10}); and
\end{example}
\begin{example}\label{ex:z213}
%{\sc Example 2.3.}
\begin{multline*}
   \z(2,\{1,3\}^n) = 4^{-n}\sum_{k=0}^n(-1)^k\z(\{4\}^{n-k})
   \bigg\{(4k+1)\z(4k+2)\\
    -4\sum_{j=1}^k\z(4j-1)\z(4k-4j+3)\bigg\},
\end{multline*}
conjectured in~\cite{BBB} and proved by Bowman and
Bradley~\cite{BowBrad}.
\end{example}

%{\sc Example 2.4.}
\begin{example}
An intriguing two-parameter, arbitrary depth evaluation involving
alternations, conjectured in~\cite{BBB} and proved herein (see
Section~\ref{S8}), is
\begin{equation}
\begin{split}
   \mu(\{-1\}^m,1,\{-1\}^n) &=
   (-1)^{m+1}\sum_{k=0}^m{n+k\choose n}A_{k+n+1}P_{m-k}\\
   &+ (-1)^{n+1}\sum_{k=0}^n{m+k\choose m}Z_{k+m+1}P_{n-k},
   \label{BBB68}
\end{split}
\end{equation}
where
\begin{equation}
   A_r := {\Li}_r(\tfrac12) = \d(r) = \sum_{k=1}^\infty
           \frac1{2^k k^r},\quad
   P_r :=  \frac{(\log 2)^r}{r!},
   \quad   Z_r := (-1)^r\z(r).
\label{A-P-Z-defn}
\end{equation}
The formula~(\ref{BBB68}) is valid for all nonnegative integers $m$
and $n$ if the divergent $m=0$ case is interpreted appropriately.
\end{example}

\begin{example}
%{\sc Example 2.5.}
If the $s_j$ are all nonpositive integers, then \[
   \bigg(\sum_{i=j}^k\nu_i\bigg)^{-s_j}
    = D_j\, \exp\bigg(-u_j\sum_{i=j}^k\nu_i\bigg),
   \quad D_j := \(-\frac{d}{du_j}\)^{-s_j}\from{u_j=0}{}.
\]
Consequently,
\begin{equation}
\begin{split}
   \l{s_1,\dots,s_k\choose b_1,\dots,b_k}
   &=\sum_{\nu_1,\ldots,\nu_k=1}^\infty \;\prod_{j=1}^k b_j^{-\nu_j}
      D_j\exp\bigg(-u_j\sum_{i=j}^k\nu_i\bigg)\\
   &=\prod_{j=1}^k D_j\sum_{\nu_j=1}^\infty b_j^{-\nu_j}
      \exp\bigg(-\nu_j\sum_{i=1}^ju_i\bigg)\\
   &=\prod_{j=1}^k D_j\bigg\{\frac{1}{b_j\exp\(\sum_{i=1}^j u_i\)-1}\bigg\}.
   \label{negative}
\end{split}
\end{equation}
In particular,~(\ref{negative}) implies
\begin{equation}
   \l{0,\dots,0\choose b_1,\dots,b_k}
   =\prod_{j=1}^k\frac1{b_j-1}.\label{zero}
\end{equation}
Despite its utter simplicity,~(\ref{zero}) points the way to
deeper waters.  For example, if we put $b_j=q^{-j}$ for each $j=1,2,\dots,k$
and note that
\[
   \l{0,0,\dots,0\choose q^{-1},q^{-2},\dots,q^{-k}}
   = \sum_{n_1>n_2>\cdots>n_k>0}\; \prod_{j=1}^k q^{n_j},
   \quad k>0,
\]
then~(\ref{zero}) implies the generating function equality
\[
   \sum_{k=0}^\infty z^k\l{0,0,\dots,0\choose q^{-1},q^{-2},\dots,q^{-k}}
   =\prod_{n=1}^\infty\(1+zq^n\)
   =\sum_{k=0}^\infty z^k\prod_{j=1}^k\frac{q^j}{1-q^j},
\]
which experts in the field of basic hypergeometric series will
recognize as a~$q$-analogue of the exponential function and a
special case of the $q$-binomial theorem, usually expressed in
the more familiar form~\cite{Gasp} as
\[
   (-zq;q)_\infty = \sum_{k=0}^\infty \frac{q^{k(k+1)/2}}{(q;q)_k} z^k.
\]
The case $k=1$, $b_1=2$, $s_1=-n$ of~(\ref{negative}) yields the
numbers~\cite{Sloane} (A000629)
\begin{equation}
   \d(-n) = \l_2(-n) = \sum_{k=1}^\infty \frac{k^n}{2^k}
   ={\Li}_{-n}(\tfrac12),
   \qquad 0\le n\in\Z,
   \label{lock}
\end{equation}
which enumerate~\cite{Kon}
the combinations of a~simplex lock having $n$ buttons,
and which satisfy the recurrence
\[
   \d(-n) = 1+\sum_{j=0}^{n-1}{n\choose j}\d(-j),\qquad 1\le n\in\Z.
\]
Also, from the exponential generating function
\[
   \sum_{n=0}^\infty \d(-n)\frac{x^n}{n!}
   = \frac{e^x}{2-e^x} = \frac{2}{2-e^x}-1,
\]
we infer~\cite{Goulson,Stanley}
that for $n\ge1$, $\tfrac12 \d(-n)$ also counts
\begin{itemize}
\item the number of ways of writing a~
sum on $n$ indices;
\item the number of
functions $f:\{1,2,\dots,n\}\to\{1,2,\dots,n\}$
such that if $j$ is in the range of $f$, then so is each value less than
or equal to $j$;
\item  the number of asymmetric generalized
weak orders on $\{1,2,\dots,n\}$;
\item the number of ordered partitions
(preferential arrangements) of $\{1,2,\dots,n\}$.
\end{itemize}

The numbers $\tfrac12 \d(-n)$ also arise~\cite{Dil}
in connection with certain constants related to the
Laurent coefficients of the Riemann zeta function.
See~\cite{Sloane} (A000670) for additional references.
\end{example}

\section{Reductions}
\label{S3}

Given the multiple polylogarithm~(\ref{polydef}), we define the
{\it depth} to be $k$, and the {\it weight} to be
$s:=s_1+\cdots+s_k$. We would like to know which sums can be
expressed in terms of lower depth sums.  When a~sum can be so
expressed, we say it {\em reduces.} Especially interesting are the
sums which completely reduce, i.e.\ can be expressed in terms of
depth-$1$ sums.  We say such sums {\em evaluate.}  The concept of
weight is significant, as all our reductions preserve it.  More
specifically, we'll see that all our reductions take the form of
a~polynomial expression which is homogeneous with respect to
weight.

There are certain sums which evidently cannot be expressed
(polynomially) in terms of lower depth sums.  Such sums are called
``irreducible''. Proving irreducibility is currently beyond the
reach of number theory. For example, proving the irrationality of
expressions like $\z(5,3)/\z(5)\z(3)$ or $\z(5)/\z(2)\z(3)$ seems
to be impossible with current techniques.

\subsection{Examples of Reductions at Specific Depths}

The functional equation (an example of
a ``stuffle'' -- see Sections~\ref{S51} through \ref{S52})
\[
   \z(s)\z(t) = \z(s,t)+\z(t,s)+\z(s+t)
\]
reduces $\z(s,s)$.

One of us (Broadhurst), using high-precision arithmetic and
integer relations finding algorithms, has found many conjectured
reductions.  One example is
\begin{equation}
\label{z413-red}
   \z(4,1,3) = - \z(5,3) + \tfrac{71}{36}\z(8)
           - \tfrac52\z(5)\z(3)+\tfrac12\z(3)^2\z(2),
\end{equation}
which expresses a multiple zeta value of depth three and weight
eight in terms of lower depth MZVs, and which was subsequently
proved. Observe that the combined weight of each term in the
reduction~(\ref{z413-red}) is preserved. The easiest proof
of~(\ref{z413-red}) uses Minh and Petitot's basis of order
eight~\cite{MinPet1}.

Broadhurst also noted that although $\z(4,2,4,2)$ is apparently irreducible
in terms of lower depth MZVs, we have the conjectured\footnote{Both
sides of~(\ref{MZVreduceEuler}) agree to at least 7900 significant figures.}
weight-12 reduction
\begin{equation}
\begin{split}
   \z(4,2,4,2)
&\eu-\tfrac{1024}{27}\l(9-,3)-\tfrac{267991}{5528}\z(12)
   -\tfrac{1040}{27}\z(9,3) -\tfrac{76}{3}\z(9)\z(3)\\
 &-\tfrac{160}9\z(7)\z(5) +2\z(6)\z(3)^2 +14\z(5,3)\z(4)\\
 &+70\z(5)\z(4)\z(3) -\tfrac16\z(3)^4
\label{MZVreduceEuler}
\end{split}
\end{equation}
in terms of lower depth MZVs {\it and} the alternating Euler sum $\l(9-,3)$.
Thus, alternating Euler sums enter quite naturally into the analysis.
And once the alternating sums are admitted, we shall see
that more general
polylogarithmic sums are required.

We remark that the depth-two sums in~(\ref{MZVreduceEuler}),
namely $\l(9-,3)$, $\z(9,3)$, and $\z(5,3)$, are almost certainly
irreducible. For example, if there are integers $c_1$, $c_2$,
$c_3$, $c_4$ (not all equal to 0) such that
$c_1\z(5,3)+c_2\pi^8+c_3\z(3)^2\z(2)+c_4\z(5)\z(3)=0$, then the
Euclidean norm of the vector $(c_1,c_2,c_3,c_4)$ is greater than
$10^{50}$.  This result can be proved computationally in a~mere
0.2 seconds on a DEC Alpha workstation using D.~Bailey's fast
implementation of the integer relation algorithm PSLQ \cite{FBA},
once we know the four input values at the precision of 200 decimal
digits. Such evaluation poses no obstacle to our fast method of
evaluating polylogs using the H\"older convolution (see
Section~\ref{S7}).

\subsection{An Arbitrary Depth Reduction}
In contrast to the specific
numerical results provided by~(\ref{z413-red}) and~(\ref{MZVreduceEuler}),
reducibility results for arbitrary sets of arguments can be obtained
if one is prepared to consider certain specific combinations of
MZVs.  The following result is typical in this respect.  It states
that, depending on the parity of the depth, either the sum or the
difference of an MZV with its reversed-string counterpart
always reduces.  Additional reductions, such as those alluded to
in Sections 1 and 2, must await the development of the theory
provided in Sections 4--7.

\begin{theorem}\label{RR}
Let $k$ be a positive integer and let $s_1,s_2,\ldots,s_k$ be
positive integers with $s_1$ and $s_k$ greater than 1.  Then the
expression
\[
   \z(s_1,s_2,\ldots,s_k)+(-1)^k\z(s_k,\ldots,s_2,s_1)
\]
reduces to lower depth MZVs.
\end{theorem}
\begin{remark} The condition on $s_1$ and $s_k$ is imposed
only to ensure convergence of the requisite sums.
\end{remark}
\begin{proof}
%To fix ideas, suppose $k=3$.
%Let $N:=\({\Z^+}\)^3 = {\Z^+}\times{\Z^+}\times{\Z^+}$
%denote the Cartesian product of three copies of the
%positive integers.  Define an additive weight-function $w:2^N\to \R$ by
%\[
%  w(A) := \sum_{(n_1,n_2,n_3)\in A} n_1^{-s_1}n_2^{-s_2}n_3^{-s_3}.
%\]
%Now if
%\begin{eqnarray*}
%   P_1 &:=&\{(n_1,n_2,n_3)\in N: n_1\le n_2\},\\
%   P_2 &:=& \{(n_1,n_2,n_3)\in N: n_2\le n_3\},
%\end{eqnarray*}
%then
%\[
%   P_1\cap P_2 = \{(n_1,n_2,n_3)\in N: n_1\le n_2\le n_3\}
%\]
%and
%\[
%   \(N\setminus P_1\)\cap\(N\setminus P_2\)
%   = \{(n_1,n_2,n_3)\in N: n_1>n_2>n_3\}.
%\]
%
%By the Inclusion-Exclusion Principle,
%\[
%   w(\(N\setminus P_1\)\cap \(N\setminus P_2\))
%   = w(N)-w(P_1)-w(P_2) +w(P_1\cap P_2),
%\]
%which is to say that
%\begin{eqnarray*}
%   \z(s_1,s_2,s_3) &=& \z(s_1)\z(s_2)\z(s_3)
%   -\z(s_3)[\z(s_2,s_1)+\z(s_2+s_1)]\\
%   &-& \z(s_1)[\z(s_3,s_2)+\z(s_3+s_2)] + \z(s_3+s_2,s_1)\\
%   &+& \z(s_3,s_2+s_1) + \z(s_3+s_2+s_1) + \z(s_3,s_2,s_1),
%\end{eqnarray*}
%i.e. $\z(s_1,s_2,s_3)-\z(s_3,s_2,s_1)$ reduces.
%
%In general, let $k$ be a positive integer and let
Let $N:=\(\Z^+\)^k$ denote the Cartesian product of $k$ copies of
the positive integers.
Define an additive weight-function $w : 2^N \to \R$ by
\[
   w(A) := \sum_{\vec n\in A}\; \prod_{j=1}^k n_j^{-s_j},
\] where the sum is over all $\vec n = (n_1,n_2,\ldots,n_k)\in
A\subseteq N.$ For each $1\le j\le k-1$, define the subset $P_j$
of $N$ by \[
   P_j := \{\vec n\in N: n_j\le n_{j+1}\}.
\]
The Inclusion-Exclusion Principle states that
\begin{equation}
   w\bigg(\bigcap_{j=1}^{k-1} N\setminus P_j\bigg)
   = \sum_{T\subseteq\{1,2,\ldots,k-1\}} (-1)^{|T|}\,
   w\bigg(\bigcap_{j\in T} P_j\bigg).
   \label{IEP}
\end{equation}
We remark that the term on the right-hand side of~(\ref{IEP})
arising from the subset $T=\{\}$ is $\z(s_1)\z(s_2)\cdots\z(s_k)$
by the usual convention for intersection over an empty set. Next,
note that the left-hand side of~(\ref{IEP}) is simply
$\z(s_1,s_2,\ldots,s_k)$.  Finally, observe that all terms on the
right-hand side of~(\ref{IEP}) have depth strictly less than
$k$---except when $T=\{1,2,\ldots,k-1\}$, which gives
\[
   (-1)^{k-1}\sum_{n_1\le n_2\le\cdots\le n_k} \;\;
      \prod_{j=1}^k n_j^{-s_j}
   = (-1)^{k-1} \z(s_k,\ldots,s_2,s_1)
        \;+\;{\mathrm{lower\; depth\; MZVs}}.
\] This latter observation completes the proof of
Theorem~\ref{RR}.
%\eop
\end{proof}

\section{Integral Representations}
\label{S4}

Writing the definition of the gamma function~\cite{Niels} in the form
\[
   r^{-s}\G(s) = \int_{1}^\infty (\log x)^{s-1} x^{-r-1}\,dx,
   \quad r>0,\quad s>0,
\]
it follows that if each $s_j>0$ and each $|b_j|\ge1$, then
\begin{equation}
\begin{split}
   \l{s_1,\dots,s_k\choose b_1,\dots,b_k} &=
   \sum_{\nu_1,\ldots,\nu_k=1}^\infty \;\prod_{j=1}^k b_j^{-\nu_j}
   \bigg(\sum_{i=j}^k\nu_i\bigg)^{-s_j}\\
   &= \sum_{\nu_1=1}^\infty
   \int_1^\infty \frac{(\log x)^{s_1-1}\,dx}
                      {\G(s_1)b_1^{\nu_1} x^{\nu_1+1}}
   \sum_{\nu_2,\ldots,\nu_k=1}^\infty \;\prod_{j=2}^k  (x b_j)^{-\nu_j}
   \bigg(\sum_{i=j}^k \nu_i\bigg)^{-s_j}\\
   &=\frac1{\G(s_1)}\int_1^\infty\frac{(\log x)^{s_1-1}}{x b_1-1}
   \l{s_2,\dots,s_k\choose x b_2,\dots,x b_k}\frac{dx}{x},
\label{recurint}
\end{split}
\end{equation}
a representation vaguely
remindful of the integral recurrence for the polylogarithm.
Repeated application of~(\ref{recurint}) yields
the $k$-dimensional integral representation
\begin{equation}
   \l{s_1,\ldots,s_k\choose b_1,\ldots,b_k}
=  \int_1^\infty\cdots\int_1^\infty\prod_{j=1}^k
    \frac{(\log x_j)^{s_j-1}dx_j}
   {\G(s_j)\big(b_j\prod_{i=1}^j x_i-1\big)x_j},
   \label{cranint1}
\end{equation}
which generalizes Crandall's integral~\cite{Crand1} for $\z(s_1,\dots,s_k)$.
An equivalent formulation of~(\ref{cranint1}) is
\begin{equation}
   \l{s_1,\ldots,s_k\choose b_1,\ldots,b_k}
  =\int_0^\infty\cdots\int_0^\infty\prod_{j=1}^k
   \frac{u_j^{s_j-1}\,du_j}
   {\G(s_j)(b_j \exp\big(\sum_{i=1}^j u_i\big) - 1)},
   \label{cranint2}
\end{equation}
the integral transforms in~(\ref{cranint2}) replacing the derivatives
in~(\ref{negative}).

Although {\em depth-dimensional integrals} such as~(\ref{cranint1})
and~(\ref{cranint2}) are attractive, they are not particularly
useful.  As mentioned previously, we are interested in reducing
the depth whenever this is possible.  However, since the weight
is an invariant of all known reductions, we seek integral representations
which respect weight invariance.  As we next show, this can be accomplished by
selectively removing logarithms from the integrand
of~(\ref{cranint1}), at the expense of increasing the number of
integrations.  At the extreme, the representation~(\ref{cranint1})
is replaced by a~{\em weight-dimensional integral} of a~rational function.

\subsection{The Partition Integral}
\label{S41} We begin with the parameters in~(\ref{polydef}).  Let
$R_1$, $R_2,\dots,R_n$ be a~(disjoint) set partition of
$\{1,2,\dots,k\}$.  Put
\[
   r_m := \sum_{i\in R_m}s_i,\quad 1\le m\le n.
\]
If $d_1,d_2,\dots,d_n$ are real numbers satisfying
$|d_m|\ge 1$ for all $m$ and $r_1d_1\ne 1$, then
\begin{eqnarray*}
   \l{r_1,\dots,r_n\choose d_1,\dots,d_n} &=&
   \sum_{\nu_1,\ldots,\nu_n=1}^\infty \;\prod_{m=1}^n
      d_m^{-\nu_m}\bigg(\sum_{j=m}^n\nu_j\bigg)^{-r_m}\\
   &=&\sum_{\nu_1,\ldots,\nu_n=1}^\infty \;\prod_{m=1}^n
    d_m^{-\nu_m}\prod_{i\in R_m}
   \bigg(\sum_{j=m}^n\nu_j\bigg)^{-s_i}\\
   &=&\sum_{\nu_1,\ldots,\nu_n=1}^\infty \;\prod_{m=1}^n
    d_m^{-\nu_m}\int_1^\infty\cdots\int_1^\infty\prod_{i\in R_m}
     \frac{\(\log x_i\)^{s_i-1}dx_i}
     {\G(s_i)x_i^{1+\nu_m+\cdots+\nu_n}}.
\end{eqnarray*}
Now collect bases with like exponents and note that
``$   \prod_{m=1}^n\prod_{i\in R_m} = \prod_{j=1}^k.$''
It follows that
\begin{equation}
\begin{split}
   \l{r_1,\dots,r_n\choose d_1,\dots,d_n}
   &=\int_1^\infty\cdots\int_1^\infty
      \bigg\{\sum_{\nu_1,\ldots,\nu_n=1}^\infty \;\prod_{m=1}^n d_m^{-\nu_m}
      \prod_{j=1}^m\prod_{i\in R_j}x_i^{-\nu_m}\bigg\}\\
      &\qquad\times
      \prod_{j=1}^k\frac{\(\log x_j\)^{s_j-1}dx_j}{\G(s_j)\,x_j}\\
   &=\int_1^\infty\cdots\int_1^\infty
      \bigg\{\prod_{m=1}^n\bigg(d_m\prod_{j=1}^m
      \prod_{i\in R_j}x_i -1\bigg)^{-1}\bigg\}\\
      &\qquad\times
      \prod_{j=1}^k
      \frac{\(\log x_j\)^{s_j-1}dx_j}{\G(s_j)\,x_j},
\label{prtnint}
\end{split}
\end{equation}
on summing the $n$ geometric series.

%{\sc Example 4.1.1.}
\begin{example}\label{ex:depth}
Taking $n=k$, we have $R_m=\{m\}$, and $r_m=s_m$ for all $1\le m\le n.$
In this case,~(\ref{prtnint}) reduces to
the depth-dimensional integral representation~(\ref{cranint1}).
\end{example}

\begin{example}
%{\sc Example 4.1.2.}
Taking $n=1$, we have $R_1=\{1,2,\dots,k\}$ and $r_1=s=\sum_{j=1}^k s_i$.
If we also put $d:=\prod_{j=1}^k d_j$, then~(\ref{prtnint})
yields the seemingly wasteful $k$-dimensional integral
\[
   \l{s\choose d} = \l{\sum_{j=1}^k s_j\choose \prod_{j=1}^k d_j}
   = \int_1^\infty\cdots\int_1^\infty
      \bigg(\prod_{j=1}^k d_j x_j-1\bigg)^{-1}
      \prod_{j=1}^k
     \frac{\(\log x_j\)^{s_j-1}dx_j}{\G(s_j)\,x_j}
\]
for a~polylogarithm of depth one.
\end{example}

%{\sc Example 4.1.3.}
\begin{example}\label{ex:weight}
Let $s_j=1$ for each $1\le j\le k$, $r_0=0$ and let
$r_1,r_2,\ldots,r_n$ be arbitrary positive integers with
$\sum_{m=1}^n r_m =k$. For $1\le m\le n$ define
\[
   R_m := \bigcup_{j=1}^{r_m} \bigg\{j+\sum_{i=1}^{m-1}r_i\bigg\}.
\]
In this case,~(\ref{prtnint}) yields a~weight-dimensional integral
of a~rational function in $k$ variables:
\begin{equation}
   \l{r_1,\dots,r_n\choose d_1,\dots,d_n}
 = \int_1^\infty\cdots\int_1^\infty
    \bigg\{\prod_{m=1}^n
    \bigg(d_m\prod_{i=1}^{u_m}x_i-1\bigg)^{-1}\bigg\}
   \prod_{j=1}^{u_n}\frac{dx_j}{x_j},
\label{ratprtn}
\end{equation}
where $u_m=\sum_{i=1}^m r_i$.  An interesting specialization
of~(\ref{ratprtn}) is
\[
\begin{split}
   \z(2,1) &= \int_1^\infty\int_1^\infty\int_1^\infty
             \frac{dx\,dy\,dz}{xyz(xy-1)(xyz-1)}
      = \int_1^\infty\int_1^\infty\int_1^\infty
         \frac{dx\,dy\,dz}{xyz(xyz-1)}\\
      &= \z(3).
\end{split}
\]
\end{example}

Although it may seem wasteful, as in Example~\ref{ex:depth} above,
to use more integrations than are required, nevertheless such
a~technique allows an easy comparison of multiple polylogarithms
having a common weight but possessing widely differing depths. For
example, from the four equations
\begin{equation}
\begin{split}
   \l{s+t\choose ab}
   &=\frac1{\G(s)\G(t)}\int_1^\infty\int_1^\infty \frac{(\log x)^{s-1}
   (\log y)^{t-1}\,dx\,dy}{(abxy-1)xy},\\
   \l{s,t\choose a,ab}
   &=\frac1{\G(s)\G(t)}\int_1^\infty\int_1^\infty \frac{(\log x)^{s-1}
   (\log y)^{t-1}\,dx\,dy}{(ax-1)(abxy-1)xy},\\
   \l{t,s\choose b,ab}
   &=\frac1{\G(s)\G(t)}\int_1^\infty\int_1^\infty \frac{(\log x)^{s-1}
   (\log y)^{t-1}\,dx\,dy}{(by-1)(abxy-1)xy},\\
   \l{s\choose a}\l{t\choose b}
   &=\frac1{\G(s)\G(t)}\int_1^\infty\int_1^\infty \frac{(\log x)^{s-1}
   (\log y)^{t-1}\,dx\,dy}{(ax-1)(by-1)xy},
\label{foureqns}
\end{split}
\end{equation}
and the rational function identity
\begin{equation}
   \frac1{(ax-1)(by-1)}
  =\frac1{abxy-1}\(\frac1{ax-1}+\frac1{by-1}+1\),
  \label{ratident}
\end{equation}
the ``stuffle'' identity (see Section~\ref{S51})
\begin{equation}
   \l{s\choose a}\l{t\choose b}
  =\l{s,t\choose a,ab} +\l{t,s\choose b,ab} +\l{s+t\choose ab}
  \label{stufident}
\end{equation}
follows immediately.  The connection between ``stuffle'' identities
and rational functions will be explained and explored more fully
in Section~\ref{S52}.

\subsection{The Iterated Integral}
\label{S42} A second approach to removing the logarithms from the
depth-dimensional integral representation~(\ref{cranint1}) yields
a~weight-dimensional iterated integral. The advantage here is that
the rational function comprising the integrand is particularly
simple.

We use the notation of Kassel~\cite{Kass} for iterated integrals.
For $j=1,2,\dots,n$, let $f_j:[a,c]\to\R$ and $\W_j := f_j(y_j)\,dy_j$.
Then
\begin{eqnarray*}
   \int_a^c \W_1\W_2\cdots\W_n
   &:=& \prod_{j=1}^n\int_a^{y_{j-1}}f_j(y_j)\,dy_j,\quad y_0:=c\\
   &=& \left\{ \begin{array}{ll}
   \int_a^c f_1(y_1) \int_a^{y_1}\W_2\cdots\W_n \,dy_1 &\mbox{if $n>0$}\\
   1 &\mbox{if $n=0$.}\\
   \end{array}\right.
\end{eqnarray*}

For each real number $b$, define a~differential $1$-form
\[
   \w_b:= \w(b) := \frac{dx}{x-b}.
\]
With this definition, the change of variable $y\mapsto1-y$
generates an involution $\w(b)\mapsto\w(1-b)$.
By repeated application of the self-evident representation
\[
   b^m m^{-s} = \int_0^b \w_0^{s-1} y^{m-1}\,dy, \qquad 1\le m\in\Z
\]
one derives from~(\ref{polydef}) that
\begin{eqnarray}
   \l{s_1,\dots,s_k \choose b_1,\dots,b_k} &=&
   \sum_{\nu_1,\ldots,\nu_k=1}^\infty \prod_{j=1}^k
      b_j^{-\nu_j}\int_0^{y_{j-1}}\w_0^{s_j-1}
      y_j^{\nu_j-1}\,dy_j,\quad y_0 := 1\nonumber\\
   &=& \prod_{j=1}^k\int_0^{y_{j-1}}\w_0^{s_j-1}
      \frac{b_j^{-1}\,dy_j}{1-b_j^{-1}y_j}\nonumber\\
   &=& (-1)^k\iu\prod_{j=1}^k\w_0^{s_j-1}\w(b_j).
   \label{iterint}
\end{eqnarray}
Letting $s:=s_1+s_2+\cdots+s_k$ denote the weight, one observes
that the representation~(\ref{iterint}) is an $s$-dimensional
iterated integral over the simplex $1>y_1>y_2>\cdots>y_s>0$.
Scaling by $q$ at each level yields the following version of the
linear change of variable formula for iterated integrals:
\begin{equation}
   \l_q{s_1,\dots,s_k\choose b_1,\dots,b_k}
   :=\l{s_1,\dots,s_k\choose qb_1,\dots,qb_k}
   =(-1)^k\int_0^{1/q}\prod_{j=1}^k\w_0^{s_j-1}\w(b_j)
   \label{scale}
\end{equation}
for any real number $q\ne0$.

Having seen that every multiple polylogarithm can be
represented~(\ref{iterint}) by a weight-dimensional iterated
integral, it is natural to ask whether the converse holds.  In
fact, any convergent iterated integral of the form
\begin{equation}
   \iu \prod_{r=1}^s \w_{\a(r)}\label{itertozeta0}
\end{equation}
can always (by collecting adjacent $\w_0$ factors -- note that
for convergence, $\a(s)\ne0$) be written in the form
\begin{equation}
   \iu \prod_{j=1}^k \w_0^{s_j-1}\w(b_j)
 = (-1)^k\l{s_1,\dots,s_k\choose b_1,\dots,b_k},
   \label{itertozeta1}
\end{equation}
where
\begin{equation}
   0\ne b_j = \a\bigg(\sum_{i=1}^j s_i\bigg).
   \label{itertozeta2}
\end{equation}

We remark that the iterated integral
representation~(\ref{iterint}) and the weight-dimen\-sional
non-iterated integral representation~(\ref{ratprtn}) of
Example~\ref{ex:weight} are equivalent under the change of
variable $x_j=y_{j-1}/y_j$, $y_0:=1$, $j=1,2,\dots,s$. In fact,
every integral representation of Section~\ref{S41} has
a~corresponding iterated integral representation under the
aforementioned transformation. For example, the depth-dimensional
integral~(\ref{cranint1}) becomes \[
   \l{s_1,\dots,s_k\choose b_1,\dots,b_k}
   = \prod_{j=1}^k \int_0^{y_{j-1}}\frac{\(\log(y_{j-1}/y_j)\)^{s_j-1}\,dy_j}
     {\G(s_j)(b_j-y_j)}.
\]
The explicit observation that MZVs are values
of iterated integrals is apparently
due to Maxim Kontsevich~\cite{Zag}.
Less formally, such representations go as far back as Euler.

\section{Shuffles and Stuffles}
\label{S5}

Although it is natural to study multiple polylogarithmic sums as
analytic objects, a~good deal can be learned from the
combinatorics of how they behave with respect to their argument
strings.

\subsection{The Stuffle Algebra}
\label{S51}

Given two argument strings $\vec s=(s_1,\dots,s_k)$ and
$\vec t=(t_1,\dots,t_r)$, we define the set
$\mbox{stuffle}(\vec s,\vec t)$ as the smallest set of
strings over the alphabet
\[
   \{s_1,\dots,s_k,t_1,\dots,t_r,\mbox{``+'',\ ``,'',\ ``('',\ ``)''}\}
\]
satisfying
\begin{itemize}
\item $(s_1,\ldots,s_k,t_1,\ldots,t_r)
         \in \mbox{stuffle}(\vec s,\vec t)$.
\item If a~string of the form $(U,s_n,t_m,V)$ is in
     stuffle$(\vec s,\vec t)$, then so are the strings
     $(U,t_m,s_n,V)$ and $(U,s_n+t_m,V)$.
\end{itemize}
Let $\vec a=(a_1,\dots,a_k)$ and
$\vec b=(b_1,\dots,b_r)$ be two strings of the same length
as $\vec s$ and $\vec t$, respectively.
We now define
\begin{equation}
   ST := ST{\vec s,\vec t\choose \vec a,\vec b}
   \label{STdef}
\end{equation}
to be the set of all pairs ${\vec u\choose \vec c}$ with
$\vec u\in\mbox{stuffle}(\vec s,\vec t)$ and $\vec c=(c_1,c_2,\dots,c_h)$
defined as follows:
\begin{itemize}
\item $h$ is the number of components of $\vec u$,

\item $c_0:= a_0 := b_0 := 1,$

\item for $1\le j\le h$, if $c_{j-1}=a_{n-1}b_{m-1}$, then
\[
   c_j := \left\{\begin{array}{lll}
         a_nb_m, &\mbox{if $u_j=s_n+t_m$},\\
         a_nb_{m-1}, &\mbox{if $u_j=s_n$},\\
     a_{n-1}b_m, &\mbox{if $u_j=t_m$}.\end{array}\right.
\]
\end{itemize}

\subsection{Stuffle Identities}
\label{S51-half}
A class of identities which we call {\em ``depth-length shuffles''} or
{\em ``stuffle identities''} is generated by a~formula for the product
of two $\l$-functions.  Consider
\[
   \l{\vec s\choose \vec a} \l{\vec t\choose \vec b}
   =\bigg\{\sum_{\nu_1,\ldots,\nu_k=1}^\infty \prod_{j=1}^k a_j^{-\nu_j}
   \bigg(\sum_{i=j}^k\nu_i\bigg)^{-s_j}\bigg\}
   \bigg\{\sum_{\xi_1,\ldots,\xi_r=1}^\infty \prod_{j=1}^r
   b_j^{-\xi_j}\bigg(\sum_{i=j}^r\xi_i\bigg)^{-t_j}\bigg\}.
\]
If we put
\begin{eqnarray*}
   &n_j := \sum_{i=j}^k\nu_i,\quad m_j := \sum_{i=j}^r\xi_i,&\\
   &a_j := \prod_{i=1}^j x_i,\quad b_j := \prod_{i=1}^j y_i,&
\end{eqnarray*}
then it follows that
\[
    \l{\vec s\choose \vec a} \l{\vec t\choose \vec b}
   =\sum_{\SMALL\begin{array}{c}
   n_1>\cdots>n_k>0\\m_1>\cdots>m_r>0\end{array}}
   \bigg(\prod_{j=1}^k x_j^{-n_j} n_j^{-s_j}\bigg)
   \bigg(\prod_{j=1}^r y_j^{-m_j} m_j^{-t_j}\bigg).
\]
Rewriting the previous expression
in terms of $\l$-functions yields the stuffle formula
\begin{equation}
   \l{\vec s\choose \vec a}\l{\vec t\choose\vec b}
   =\sum\l{\vec u\choose \vec c},
   \label{stuffle}
\end{equation}
where the sum is over all pairs of strings
${\vec u\choose \vec c}\in ST{\vec s,\vec t\choose \vec a,\vec b}$.

%{\sc Example 5.2.1.}
\begin{example}
\[
\begin{split}
   \l{r,s\choose a,b}\l{t\choose c}
   &=\l{r,s,t\choose a,b,bc}+\l{r,s+t\choose a,bc}\\
   & +\l{r,t,s\choose a,ac,bc}
   +\l{r+t,s\choose ac,bc}+\l{t,r,s\choose c,ac,bc}.
\end{split}
\]
When specialized to MZVs, this example produces the identity
\[
   \z(r,s)\z(t) = \z(r,s,t)+\z(r,s+t)+\z(r,t,s)+\z(r+t,s)+\z(t,r,s).
\]
\end{example}

The term ``stuffle'' derives from the manner in which
the two (upper) strings
are combined.  The relative order of the two strings is preserved
(shuffles), but elements of the two strings may also be shoved
together into a~common slot (stuffing), thereby reducing the depth.

\subsection{Stuffles and Partition Integrals}
\label{S52}

In Section~\ref{S41}, an example was given in which a~stuffle
identity~(\ref{stufident}) was seen to arise from a~corresponding
rational function identity~(\ref{ratident}) and certain partition
integral representations~(\ref{foureqns}).  This is by no means
an isolated phenomenon.  In fact, we shall show that {\it every}
stuffle identity is a~consequence of the partition
integral~(\ref{prtnint})
applied to a~corresponding rational function identity.

\begin{theorem}
\label{T1}
Every stuffle identity is equivalent to a~rational
function identity, via the partition integral.
\end{theorem}

Before proving Theorem~\ref{T1}, we define a class of rational
functions, and prove they satisfy a~certain rational function
identity. Let $\vec s=(s_1,\dots,s_k)$ and $\vec
t=(t_1,\dots,t_r)$ be vectors of positive integers, and let
$\vec\a=(\a_1,\dots,\a_k)$ and $\vec\b=(\b_1,\dots,\b_r)$ be
vectors of real numbers.  As in~(\ref{STdef}), put
\[
   ST = ST{\vec s,\vec t\choose \vec \a,\vec\b},
\]
and define
\[
   T = T(\vec\a,\vec\b) :=
   \bigg\{\vec\g: {\vec u\choose \vec\g}\in ST\bigg\}.
\]
Let $f:T\to\Q[\g_1,\g_2,\dots]$ be defined by
\begin{equation}
   f(\g_1,\dots,\g_h) := \prod_{j=1}^h (\g_j-1)^{-1}.
   \label{fdef}
\end{equation}
Then we have the following lemma.

\begin{lemma}
\label{L1}
Let $f$ be defined as in~(\ref{fdef}).  Then
\[
   f(\vec\a)f(\vec\b) = \sum_{\vec\g\in T(\vec\a,\vec\b)}f(\vec\g).
\]
\end{lemma}

\begin{demo}{Proof of Lemma~\ref{L1}} Apply~(\ref{stuffle}) with
$\vec a=\vec\a$ and $\vec b=\vec\b$.  In view of~(\ref{zero}), the
lemma follows on taking $\vec s$ and $\vec t$ to be zero vectors
of the appropriate lengths.
%\eop
\end{demo}

\begin{demo}{Proof of Theorem~\ref{T1}} Let $\vec s$, $\vec t$, $\vec a$,
and $\vec b$ be as in~(\ref{stuffle}).  Let $\vec\a$ and $\vec\b$
be given by
\[
   \a_j := a_j \prod_{i=1}^j x_i,\quad
   \b_j := b_j \prod_{i=1}^j y_i.
\]
Applying Lemma~\ref{L1} and the partition integral
representation~(\ref{prtnint}) to the depth-dimensional
integral~(\ref{cranint1}) yields
\[
\begin{split}
   \l{\vec s\choose \vec a}\l{\vec t\choose \vec b}
   &= \bigg\{\int_1^\infty\cdots\int_1^\infty f(\vec\a)
             \prod_{j=1}^k\frac{(\log x_j)^{s_j-1}dx_j}
                          {\G(s_j)\,x_j}\bigg\}\\
   &\qquad\times
      \bigg\{\int_1^\infty\cdots\int_1^\infty f(\vec\b)
             \prod_{j=1}^r\frac{(\log y_j)^{t_j-1}dy_j}
             {\G(t_j)\,y_j}\bigg\}\\
   &= \int_1^\infty\cdots\int_1^\infty
      \sum_{\vec\g\in T(\vec\a,\vec\b)} f(\vec\g)
      \bigg\{\prod_{j=1}^k\frac{(\log x_j)^{s_j-1}dx_j}
      {\G(s_j)\,x_j}\bigg\}\\
   &\qquad\times
      \bigg\{\prod_{j=1}^r
      \frac{(\log y_j)^{t_j-1}dy_j}{\G(t_j)\,y_j}\bigg\}\\
   &= \sum_{\vec\g\in T(\vec\a,\vec\b)}\int_1^\infty\cdots\int_1^\infty
       f(\vec\g)\bigg\{\prod_{j=1}^k\frac{(\log x_j)^{s_j-1}dx_j}
     {\G(s_j)\,x_j}\bigg\}\\
   &\qquad\times
      \bigg\{\prod_{j=1}^r
      \frac{(\log y_j)^{t_j-1}dy_j}{\G(t_j)\,y_j}\bigg\}\\
   &=\sum_{{\vec u\choose \vec c}\in ST{\vec s,\vec t\choose \vec a,\vec b}}
       \l{\vec u\choose \vec c},
\end{split}
\]
as required.
%\eop
\end{demo}

\subsection{The Shuffle Algebra}
\label{S53}
As opposed to depth-length shuffles, or stuffles, which arise
from the definition~(\ref{polydef}) in terms of sums,
the iterated integral representation~(\ref{iterint})
gives rise to what are called
{em ``weight-length shuffles'',} or simply ``shuffles''.
Weight-length shuffles take the form
\begin{equation}
   \iu \W_1\W_2\cdots\W_n \iu\W_{n+1}\W_{n+2}\cdots\W_{n+m}
 = \sum \iu \W_{\sigma(1)}\W_{\sigma(2)}\cdots\W_{\sigma(n+m)},
   \label{shuff}
\end{equation}
where the sum is over all ${n+m\choose n}$ permutations $\sigma$ of the set
$\{1,2,\ldots,n+m\}$ which satisfy $\sigma^{-1}(i)<\sigma^{-1}(j)$ for
all $1\le i<j\le n$ and $n+1\le i<j\le n+m$.  In other words, the sum
is over all $(n+m)$-dimensional iterated integrals in which the relative
orders of the two strings of $1$-forms $\W_1,\dots,\W_n$ and
$\W_{n+1},\dots,\W_{n+m}$ are preserved.

%{\sc Example 5.4.1.}
\begin{example}
\begin{eqnarray*}
   \z(2,1)\z(2)
   &=& -\iu \w_0\w_1^2 \iu \w_0\w_1\\
   &=& -6\iu\w_0^2\w_1^3 - 3\iu\w_0\w_1\w_0\w_1^2
       -\iu\w_0\w_1^2\w_0\w_1\\
   &=& 6\z(3,1,1)+ 3\z(2,2,1) +\z(2,1,2).
\end{eqnarray*}
In contrast, the stuffle formula gives
\[
   \z(2,1)\z(2) = 2\z(2,2,1) + \z(4,1) +\z(2,3)+\z(2,1,2).
\]
\end{example}

Note that weight-length shuffles preserve both depth and weight.
In other words, the depth (weight) of each term which occurs in
the sum over shuffles is equal to the combined depth (weight) of
the two multiple polylogarithms comprising the product.

Though it may appear that the shuffles form a rather trivial class
of identities satisfied by iterated integrals, it is worth
mentioning that the second proof of Zagier's conjecture (see
Corollary~\ref{C2} of Section~\ref{S102}) uses little more than
the combinatorial properties of shuffles~\cite{BBBLc}.  In
addition, both shuffles and stuffles have featured in the
investigations of other authors in related
contexts~\cite{Hoff2,Hoff3,Hoff4,Min1,Min2,MinPet2,MinPet1,MinPetHoev2,MinPetHoev1,MinPetHoev3,Reut}.

\section{Duality}
\label{S6}
In~\cite{Hoff1}, Hoffman defines an involution on strings $s_1,\dots,s_k$.
The involution coincides with a~notion we refer to as duality.
The duality principle states that two MZVs coincide whenever their argument
strings are dual to each other, and (as noted by Zagier~\cite{Zag})
follows readily from the iterated integral representation.
In~\cite{DJB2}, Broadhurst generalized the notion of duality
to include relations between iterated integrals involving the sixth root
of unity; here we allow arbitrary complex values of $b_j$.
Thus, we find that the duality principle easily extends to
multiple polylogarithms, and in this more general setting, has
far-reaching implications.

\subsection{Duality for Multidimensional Polylogarithms}
\label{S61}
We begin with the iterated integral representation~(\ref{iterint}) of
Section~\ref{S42}.
Reversing the order of the omegas and replacing each integration
variable $y$ by its complement $1-y$ yields the dual iterated integral
representation
\begin{equation}
   \l{s_1,\dots,s_k\choose b_1,\dots,b_k}
   = (-1)^{s+k}\iu\prod_{j=k}^1 \w(1-b_j)\w_1^{s_j-1},
   \label{dualint}
\end{equation}
where again $s=s_1+\cdots+s_k$ is the weight.

%{\sc Example 6.1.1.}
\begin{example}
Using~(\ref{polydef}),~(\ref{iterint}), and~(\ref{dualint}),
we have
\begin{eqnarray*}
   \l{2,1\choose1,-1}
   = \int_0^1 \w(0)\,\w(1)\,\w(-1)
   =-\int_0^1\w(2)\,\w(0)\,\w(1)
   =-\l{1,2\choose2,1},
\end{eqnarray*}
which is to say that
\[
   \sum_{n=1}^\infty\frac1{n^2}\sum_{k=1}^{n-1}\frac{(-1)^k}k
   =-\sum_{n=1}^\infty\frac1{n2^n}\sum_{k=1}^{n-1}\frac{2^k}{k^2},
\]
a result that would doubtless be difficult to prove by na{\"\i}ve series
manipulations alone.
\end{example}

When $b_1=b_2=\dots=b_k=b$, the two dual iterated integral
representations~(\ref{iterint}) and~(\ref{dualint})
simplify as follows:
\begin{equation}
   \l_b(s_1,\dots,s_k) = (-1)^k\iu\prod_{j=1}^k \w_0^{s_j-1}\w(b)
   = (-1)^{s+k}\iu\prod_{j=k}^1 \w(1-b)\w_1^{s_j-1}.
   \label{simpledual}
\end{equation}
A somewhat more symmetric version of~(\ref{simpledual}) is
\begin{eqnarray}
   (-1)^m\l_b(s_1+2,\us^{r_1},\dots,s_m+2,\us^{r_m})
   &=&(-1)^r\iu\prod_{j=1}^m\w_0^{s_j+1}\w_b^{r_j+1}\nonumber\\
   &=&(-1)^s\iu\prod_{j=m}^1\w_{1-b}^{r_j+1}\w_1^{s_j+1},
   \label{symdual}
\end{eqnarray}
where $r:=\sum_j r_j$ and, as usual, $s:=\sum_j s_j$.

\subsection{Duality for Unsigned Euler Sums}
\label{S62}
Taking $b=1$ in~(\ref{symdual}), we deduce the
{\em MZV duality
formula} (cf.~\cite{Kass} p.~483)
\begin{equation}
   \z(s_1+2,\us^{r_1},\dots,s_m+2,\us^{r_m})
   = \z(r_m+2,\us^{s_m},\dots,r_1+2,\us^{s_1})\label{MZVdual}
\end{equation}
for multidimensional unsigned Euler sums, i.e.\ multiple zeta
values (MZVs) .

%{\sc Example 6.2.1.}
\begin{example}
MZV duality~(\ref{MZVdual}) gives Euler's evaluation $\z(2,1) =
\z(3),$ as well as the generalizations $\z(\{2,1\}^n) =
\z(\{3\}^n),$ and $\z(2,\us^n) = \z(n+2)$, valid for all
nonnegative integers $n$.
\end{example}

In~\cite{YOhno} a beautiful extension of MZV duality~(\ref{MZVdual})
is given, which also subsumes the so-called sum
formula
\[
   \sum_{\stackrel{\scriptstyle n_j>\delta_{j,1}}{N=\Sigma_j n_j}}
   \z(n_1,n_2,\ldots,n_k)=\z(N),
\] conjectured independently by C. Moen~\cite{Hoff1} and M.
Schmidt~\cite{Mar}, and subsequently proved by A.
Granville~\cite{AGran}. We refer the reader to Dr. Ohno's article
for details.

The duality principle has an enticing converse, namely that
{\em two MZVs
with distinct argument strings are equal only if the argument strings
are dual to each other.}
Unfortunately, although the numerical (and symbolic)
evidence in support of this
converse statement is overwhelming, it still remains to be proved.
In the case of self-dual strings, the conjectured converse
of the duality principle implies that such a MZV can equal no other
MZV; moreover we find that certain of these
completely reduce, i.e. evaluate entirely in terms of (depth-one)
Riemann zeta functions.

%{\sc Example 6.2.2.}
\begin{example}\label{ex:Zag}
The following self-dual evaluation, previously conjectured by Don
Zagier~\cite{Zag}
\[
   \z(\{3,1\}^n) = 4^{-n}\z(\{4\}^n) = \frac{2\pi^{4n}}{(4n+2)!},
   \qquad 0\le n\in\Z,
\]
is proved herein (see Section~\ref{S10}).
\end{example}

\begin{example}
The evaluation
\begin{multline*}
   \z(2,\{1,3\}^n) = 4^{-n}\sum_{k=0}^n(-1)^k\z(\{4\}^{n-k})
   \bigg\{(4k+1)\z(4k+2)\\
    -4\sum_{j=1}^k\z(4j-1)\z(4k-4j+3)\bigg\},
    \qquad 0\le n\in\Z
\end{multline*}
conjectured in~\cite{BBB} and
recently proved by Bowman and Bradley~\cite{BowBrad}
is also self-dual.
\end{example}

\begin{example}
The self-dual two-parameter
generalization of Example~\ref{ex:Zag}
\[
   \z(\{2\}^m,\{3,\{2\}^m,1,\{2\}^m\}^n)
   \eu \frac{2(m+1)\pi^{4(m+1)n+2m}}{(2(m+1)(2n+1))!},
   \qquad 0\le m,n\in\Z,
\]
remains to be proved.
\end{example}

We conclude this section with
the following result, since the special case $p=1$ has
some bearing on the MZV duality formula~(\ref{MZVdual}).

\begin{theorem}
\label{T2}
 Let $|p|\ge 1$.  The double generating function equality
\[
1 -  \sum_{m=0}^\infty\sum_{n=0}^\infty x^{m+1}y^{n+1}\l_p(m+2,\us^n)
  = {}_2F_1\(\begin{array}{c}y, -x\\ 1-x\end{array}
         \left|\phantom{\displaystyle\int}\!\!\!\!\frac1p\right.\)
\]
holds.
\end{theorem}

\begin{proof} By definition~(\ref{polysubdef}) of $\l_p$,
\begin{eqnarray*}
   \sum_{m=0}^\infty\sum_{n=0}^\infty x^{m+1}y^{n+1}\l_p(m+2,\us^n)
   &=& y\sum_{m=0}^\infty x^{m+1}\sum_{k=1}^\infty
      \frac1{k^{m+2}p^k}\prod_{j=1}^{k-1} \(1+\frac{y}j\)\\
   &=& \sum_{m=0}^\infty x^{m+1}\sum_{k=1}^\infty
      \frac{(y)_k}{k^{m+1}k!p^k}\\
   &=& \sum_{k=1}^\infty\frac{(y)_k}{k! p^k}\(\frac{x}{k-x}\)\\
   &=& -\sum_{k=1}^\infty\frac{(y)_k(-x)_k}{k! p^k(1-x)_k}\\
   &=& 1-{}_2F_1\(\begin{array}{c}y, -x\\ 1-x\end{array}
       \left|\phantom{\displaystyle\int}\!\!\!\!\frac1p\right.\)
\end{eqnarray*}
as claimed.
%\eop
\end{proof}

\begin{remarks} In~\cite{BBB} it was noted that the $p=1$ case
of Theorem~\ref{T3}
is equivalent to the $m=1$ case of MZV duality~(\ref{MZVdual})
via the invariance of
\begin{equation}
\begin{split}
   {}_2F_1\(\begin{array}{c}y, -x\\ 1-x\end{array}
      \left|\phantom{\displaystyle\int}\!\!\!\!1\right.\)
   &=\frac{\G(1-x)\G(1-y)}{\G(1-x-y)}\\
   &= \exp\bigg\{\sum_{k=2}^\infty\(x^k+y^k-(x+y)^k\)\frac{\z(k)}{k}\bigg\}
\label{hyperGauss}
\end{split}
\end{equation}
with respect to the interchange of $x$ and $y$.  However, it appears
that this observation can be traced back to Drinfeld~\cite{Drin}.  In
connection with his work on series of Lie brackets, Drinfeld encountered
a scaled version of the exponential series above, and showed that the
coefficients of the double generating function satisfy $c_{mn}=c_{nm}$
and $c_{m0}=c_{0m}$ evaluates to $\z(m+2)$, up to
a~so-called Oppenheimer factor
which we omit (\cite{Kass}, p.~468).  In our notation, this is essentially
the statement that $\z(m+2,\us^n)=\z(n+2,\us^m)$.

Note that Theorem~\ref{T3} in conjunction with~(\ref{hyperGauss})
shows that $\z(m+2,\us^n)$ completely reduces (i.e. is expressible
solely in terms of depth-1 Riemann zeta values) for all
nonnegative integers $m$ and $n$.  In particular, the coefficient
of $x^{m-1}y^2$ gives Euler's formula (Example~\ref{ex:z21});  and
taking the coefficient of $x^{m-1}y^3$ provides a much simpler
derivation of Markett's formula~\cite{Mar} for $\z(m,1,1)$, $m\ge
2$.  Thus, the complete reducibility of $\z(m+2,\us^n)$ is
a~simple consequence of the instance~(\ref{hyperGauss}) of Gauss's
${}_2F_1$ hypergeometric summation
theorem~\cite{AS,WNBailey,Slater}. Wenchang Chu~\cite{CHU} has
elaborated on this idea, applying additional hypergeometric
summation theorems to evaluate a wide variety of depth-2 sums,
including nonlinear (cf.~\cite{Flaj}) sums.

It would be interesting to know if there is a~generating function formulation
of MZV duality at full strength~(\ref{MZVdual}).  Presumably, it would
involve an analogue of Drinfeld's associator in $2m$ non-commuting
variables.
\end{remarks}

\subsection{Duality for Unit Euler Sums}
\label{S63}
Recall the $\d$-function was defined~(\ref{def-zeta-lambda})
as the nested sum extension of the polylogarithm at one-half:
\begin{equation}
   \d(s_1,\dots,s_k) :=\l{s_1,\dots,s_k\choose 2,\dots,2} =
   \sum_{\nu_1,\ldots,\nu_k=1}^\infty \;\prod_{j=1}^k 2^{-\nu_j}
   \bigg(\sum_{i=j}^k\nu_i\bigg)^{-s_j}.
   \label{deltadef}
\end{equation}
Due to its geometric rate of convergence, $\d$-values
can be computed to high precision relatively
quickly.  On the other hand, the unit Euler $\mu$-sums~(\ref{unitEuler})
converge extremely slowly when the $b_j$ all lie on the unit circle.
In particular, the slow convergence of the unit $(\pm 1)$
argument $\mu$-sums initially confounded our efforts to create a
data-base of numerical evaluations from which to form viable
conjectures.
Nevertheless, there is a~close relationship between the $\d$-sums
and the $\mu$-sums, as we shall presently see.

Taking $b=2$ in~(\ref{symdual}), we deduce the ``delta-to-unit-mu''
duality formula
\begin{multline}
   \d(s_1+2,\us^{r_1},\ldots,s_m+2,\us^{r_m})\\
   = (-1)^{r+m}
   \mu(\{-1\}^{r_m+1},\us^{s_m+1},\dots,\{-1\}^{r_1+1},\us^{s_1+1}).
   \label{deltadualsym}
\end{multline}
Thus, every convergent unit $(\pm 1)$
argument $\mu$-sum can be expressed as a~
(rapidly convergent) $\d$-sum.  The converse follows from the more
general, but less symmetric formula, arising from~(\ref{simpledual}):
\begin{equation}
   \d(s_1,\ldots,s_k)
  =(-1)^k\mu(-1,\us^{s_k-1},\dots,-1,\us^{s_1-1}).
  \label{deltadual}
\end{equation}

%{\sc Example 6.3.1.}
\begin{example}
\[
   \d(1) = \sum_{\nu=1}^\infty\frac1{\nu 2^\nu} = -\log(\tfrac12)
         =  \sum_{\nu=1}^\infty\frac{(-1)^{\nu+1}}{\nu} = -\mu(-1),
\]
and more generally, for all nonnegative integers $n$, we have
\begin{equation}
   \d(n+1) = \sum_{\nu=1}^\infty\frac1{\nu^{n+1}2^\nu}={\Li}_{n+1}(\tfrac12)
           = -\mu(-1,\us^n).\label{delta-Li}
\end{equation}
\end{example}

%{\sc Example 6.3.2.}
\begin{example}
For all nonnegative integers $n$,
\begin{eqnarray}
   \d(\us^n)
   &=& (-1)^n\mu(\{-1\}^n) = (\log 2)^n/n!,
       \label{deltaus}\\
   \d(2,\us^n) &=& (-1)^{n+1}\mu(\{-1\}^{n+1},1),
   \label{delta2us}
\end{eqnarray}
and more generally,
\[
   \d(\us^m,2,\us^n) = (-1)^{m+n+1}\mu(\{-1\}^{n+1},1,\{-1\}^m),
   \qquad 0\le m,n\in\Z.
\]
\end{example}

%{\sc Example 6.3.3.}
\begin{example}
%\begin{eqnarray*}
\[
   \d(1,n+1) = \mu(-1,\us^n,-1),\qquad 0\le n\in\Z,
   %\d(1,n) &=& \int_0^{1/2}\frac{{\Li}_n(z)}{1-z}\,dz.
   %\quad {\bf Proof\; by\; interchanging\; summation\; order}
\]
%\end{eqnarray*}
and in particular, remembering~(\ref{A-P-Z-defn},\ref{lock},\ref{delta-Li})
that $\d(r)={\Li}_r(\tfrac12)$, we have
\begin{eqnarray*}
   \d(1,0) &=& 1-\log 2 = 1-\d(1),\\
   \d(1,2) &=& \tfrac57\d(2)\d(1) - \tfrac27\d(3) + \tfrac5{21}\d^3(1).
         %\quad{\bf Proof:\; See\;\cite{Lewin1}}\\
\end{eqnarray*}
Integer relation  searches (see \cite{BL} or \cite{BBB} for details)
have failed to find a~similar formula for $\d(1,4)$.
However,
\[
   2\d(1,2n-1) = \sum_{j=1}^{2n-1} (-1)^{j+1}\d(j)\d(2n-j),
   \qquad 1\le n\in\Z.
\]
Also,
\[
   \d(1,-n) = \sum_{\nu=0}^n{n\choose \nu}\frac{B_{n-\nu}\d(-\nu)}{\nu+1},
   \qquad 1\le n\in\Z,
\]
where the $\d(-\nu)$ are the simplex lock numbers~(\ref{lock})
and the $B_\nu$ are the Bernoulli numbers~\cite{AS}.
More generally, if $n_1$ is a~positive integer and $n_2,n_3,\dots,n_r$
are all nonnegative integers, then
\[
   \d(s,-n_r,\dots,-n_2,-n_1)
   =\bigg\{\prod_{j=1}^r\sum_{\nu_j=0}^{\tau_j}A(\nu_j)\bigg\}
   \d(s-\nu_r-1),\quad s\in\C,
\]
where
\[
   \tau_j := n_j+\nu_{j-1}+1,\quad
   A(\nu_j) := \frac1{\nu_j+1}{\tau_j\choose \nu_j}B_{\tau_j-\nu_j},\quad
   \nu_0 :=-1.
\]
%{\tt Does $\d(a,b)+\d(b,a)$ always reduce?  I don't think so.}
\end{example}

\section{The H\"older Convolution}
\label{S7} Richard Crandall~\cite{Crand2} (see
also~\cite{CranBuh}) describes a practical method for fast
evaluation of MZVs.  Here, we develop an entirely different
approach which is based on the fact that any multiple
polylogarithm can be expressed as a convolution of rapidly
convergent multiple polylogarithms.  We have used such
representations to compute otherwise slowly convergent alternating
Euler sums and (unsigned) MZVs to precisions in the thousands of
digits.  Lest this strike the reader as perhaps an excessive
exercise in recreational computation, consider that many of our
results were discovered via exhaustive numerical
searches~\cite{BBB} for which even hundreds of digits of precision
were insufficient, depending on the type of relation
sought~\cite{BL}.

A publicly available implementation of our technique is briefly described
in Section~\ref{sec-EZ}.
There are also interesting theoretical considerations which
we have only begun to explore.
See equations~(\ref{zetaLilog})--(\ref{Li2-half}) below
for a taste of what is possible.

\subsection{Derivation and Examples}

We have seen how multiple polylogarithms with unit arguments can
be expres\-sed in terms of rapidly convergent $\d$-sums. What if
the arguments are not necessarily units? In the iterated integral
representation~(\ref{iterint}) the domain $1>y_j>y_{j+1}>0$ in
$s=\sum_j s_j$ variables splits into $s+1$ parts.  Each part is
a~product of regions $1>y_j>y_{j+1}>1/p$ for the first $r$
variables, and $1/p>y_j>y_{j+1}>0$ for the remaining $s-r$
variables. Next, $y_j\mapsto 1-y_j$ replaces an integral of the
former type by one of the latter type, with $1/p$ replaced by $1/q
:= 1-1/p$.

Motivated by these observations, we consider the string of differential
$1$-forms which occurs in the integrand of the iterated integral
representation~(\ref{iterint}) and define
\[
   \a_r := \left\{\begin{array}{ll}
             b_j, & \mbox{if $r=\sum_{i=1}^j s_i$}\\
             0,  & \mbox{otherwise.}\\
                  \end{array}\right.
\]
Then
\begin{eqnarray}
   \l{s_1,\ldots,s_k\choose b_1,\ldots,b_k}
&=&(-1)^k\iu \prod_{r=1}^s\w(\a_r)\nonumber\\
&=&\sum_{r=0}^s (-1)^{r+k}\bigg\{\int_0^{1/q}\prod_{j=r}^1 \w(1-\a_j)\bigg\}
   \bigg\{\int_0^{1/p}\prod_{j=r+1}^s \w(\a_j)\bigg\}.
   \label{holder}
\end{eqnarray}
Thus, by means of~(\ref{itertozeta0}),~(\ref{itertozeta1}),
and~(\ref{itertozeta2}), we have expressed the general multiple
polylogarithm as a~convolution of $\l_p$ with $\l_q$ for any $p$,
$q$ such that the H\"older condition $1/p+1/q=1$ is satisfied. For
this reason, we refer to~(\ref{holder}) as the {\em H\"older
convolution.}  Note that the H\"older convolution generalizes
duality~(\ref{dualint}) for multiple polylogarithms, as can be
seen by letting $p$ tend to infinity so that~(\ref{scale})
$\l_p\to0$, and $q\to1$.

{\sc MZV Example.}
For any $p>0$, $q>0$ with $1/p+1/q=1$,
\begin{eqnarray*}
   \z(2,1,2,1,1,1)
   &=& \l_p(2,1,2,1,1,1)+\l_p(1,1,2,1,1,1)\l_q(1)\\
   &&+\l_p(1,2,1,1,1)\l_q(2)+\l_p(2,1,1,1)\l_q(3)\\
   &&+\l_p(1,1,1,1)\l_q(1,3)+\l_p(1,1,1)\l_q(2,3)+\l_p(1,1)\l_q(3,3)\\
   &&+\l_p(1)\l_q(4,3)+\l_q(5,3)\\
   &=& \z(5,3).
\end{eqnarray*}

The pattern should be clear.  For $1\le j\le m$, define the concatenation
products
\begin{eqnarray*}
   \vec a_j &:=& \Cat{i=j}{m} \{s_i+2,\us^{r_i}\}
   = \{s_j+2,\us^{r_j},\dots,s_m+2,\us^{r_m}\},\\
   \vec b_j &:=& \Cat{i=j}{1} \{r_i+2,\us^{s_i}\}
   = \{r_j+2,\us^{s_j},\dots,r_1+2,\us^{s_1}\},
\end{eqnarray*}
and $\vec a_{m+1} := \vec b_0 := \{\}$.  Then the H\"older convolution
for the general MZV case is given by
\begin{eqnarray}
   \z(\vec a_m) &=&
   \sum_{j=1}^m\left\{\sum_{t=0}^{s_j+1}\l_p(s_j+2-t,\us^{r_j},\vec a_{j+1})
   \l_q(\us^t,\vec b_{j-1})\right.\nonumber\\
   &&\qquad\left.+\sum_{\nu=1}^{r_j}\l_p(\us^{\nu},\vec a_{j+1})
   \l_q(r_j+2-\nu,\us^{s_j},\vec b_{j-1})\right\}+ \l_q(\vec b_m)
   \label{MZVholder}\\
    &=& \z(\vec b_m)\nonumber.
\end{eqnarray}
Of course, $\vec a_m$ and $\vec b_m$ are the dual strings in the
MZV duality formula~(\ref{MZVdual}). Since the sums $\l_p$
converge geometrically, whereas MZV sums converge only
polynomially,~(\ref{MZVholder}) provides an excellent method of
computing general MZVs to high precision with the optimal
parameter choice $p=q=2$.  For rapid computation of general
multiple polylogarithms, it is simplest to use the H\"older
convolution~(\ref{holder}) directly, translating the iterated
integrals into geometrically convergent sums on a~case by case
basis, using~(\ref{iterint}).

{\sc Alternating Example.}
\begin{eqnarray*}
   \l(2,1-)
   &=& \iu\w(0)\w(1)\w(-1)\\
   &=& \int_0^{1/p}\w(0)\w(1)\w(-1)
   -\int_0^{1/q}\w(1)\int_0^{1/p}\w(1)\w(-1)\\
   &&+\int_0^{1/q}\w(0)\w(1)\int_0^{1/p}\w(-1)
   -\int_0^{1/q}\w(2)\w(0)\w(1)\\
   &=& \l_p(2,1-)+\l_p(1,1-)\l_q(1)+\l_p(1-)\l_q(2)
   - \l_q{1,2\choose 2,1}\\
   &=& - \l{1,2\choose 2,1}.
\end{eqnarray*}
Although we could now work out the explicit form of
the analogue to (\ref{MZVholder}) in the alternating case,
the resulting formula is too complicated in relation to
its importance to justify including here.

In addition to the impressive computational implications already
outlined, the H\"older convolution~(\ref{holder}) gives new
relationships between multiple polylogarithms, providing a~path to
understanding certain previously mysterious evaluations. For
example, taking $p=q=2$ shows that every MZV of weight $s$ can be
written as a~weight-homogeneous convolution sum involving $2s$
$\d$-functions. Furthermore, employing the weight-length shuffle
formula~(\ref{shuff}) to each product shows that every MZV of
weight $s$ is a~sum of $2^s$ (not necessarily distinct)
$\d$-values, each of weight $s$, and each appearing with unit
$(+1)$ coefficient.  In particular, this shows that the vector
space of rational linear combinations of MZVs is spanned by the
set of all $\d$-values. Thus,
\[
\begin{split}
   \z(3) &= -\int_0^{1/2}\w_0\w_0\w_1 +\int_0^{1/2}\w_1\int_0^{1/2}
   \w_0\w_1-\int_0^{1/2}\w_1\w_1\int_0^{1/2}\w_1\\
   &\qquad+\int_0^{1/2}\w_0\w_1\w_1\\
   &= \d(3)+\int_0^{1/2}(\w_1\cdot\w_0\w_1+\w_0\cdot\w_1\cdot\w_1
       +\w_0\w_1\cdot\w_1)\\
   &\qquad -\int_0^{1/2}(\w_1\w_1\cdot\w_1+\w_1\cdot\w_1\cdot\w_1
       +\w_1\cdot\w_1\w_1) +\d(2,1)\\
   &= \d(3) + \d(1,2)+\d(2,1)+\d(2,1)+\d(1,1,1)+\d(1,1,1)+\d(1,1,1)+\d(2,1).
\end{split}
\]

{\sc Polylog Example.}
Applying~(\ref{holder})
to $\z(n+2)$, with $p=q=2$ provides a~lovely closed form
for $\d(2,\us^n).$   Indeed,
\begin{equation}
   \z(n+2) = \d(2,\us^n) + \sum_{r=1}^{n+2}\d(r)\d(\us^{n+2-r}).
   \label{zetaLilog}
\end{equation}
The desired closed form follows
after rearranging the previous equation~(\ref{zetaLilog}) and
applying the definition~(\ref{deltadef}) and the result~(\ref{deltaus})
in the form
$\d(r)={\Li}_r(\tfrac12)$ and $\d(\us^r)=(\log 2)^r/r!$,
respectively.

%{\sc Example.}
\begin{example}
Putting $n=1$ in~(\ref{zetaLilog}) gives~\cite{BerndtI}
\begin{equation}
\label{z3-zetaLilog}
   \z(3) = \sum_{n=1}^\infty\frac1{n^3}
   = \frac{1}{12}\pi^2 \log(2)+\sum_{n=1}^\infty\frac1{2^n n^2}\sum_{j=1}^n\frac1j.
\end{equation}
\end{example}

In fact, formula~(\ref{zetaLilog}) is non-trivial even
when $n=0$.
Putting $n=0$ in~(\ref{zetaLilog}) gives the classical
evaluation of the dilogarithm at one-half:
\begin{equation}
   2{\Li_2}(\tfrac12) = \z(2)-(\log 2)^2\quad{\mathrm i.e.}\quad
   \sum_{n=1}^\infty\frac1{2^n n^2}=\tfrac1{12}\pi^2-\tfrac12(\log 2)^2.
   \label{Li2-half}
\end{equation}

Differentiation of~(\ref{holder}) with respect to the parameter
$p$ provides another avenue of pursuit
which has not yet been fully explored.  We have used this approach
to derive $\d(0,\us^n) = \d(\us^n),$ but in fact,
removing the initial zero is trivial from first principles.

\subsection{EZ Face}
\label{sec-EZ}

A fast program for evaluating MZVs (as well as arithmetic
expressions containing them) based on the H\"older convolution
formula~(\ref{MZVholder}) has been developed at the
CECM\footnote{Centre for Experimental and Constructive
Mathematics, Simon Fraser University.}, and is available for
public use via the World Wide Web interface called ``EZ Face'' (an
abbreviation for Euler Zetas interFace) at the URL \vskip.1in
\centerline{{\tt http://www.cecm.sfu.ca/projects/EZFace/}}
\vskip.1in
\noindent
This publicly accessible interface currently allows one to evaluate the sums
\[
{\tt z}(s_1,\ldots,s_k)
  :=\sum_{n_1>\cdots>n_k} \;
    \prod_{j=1}^k n_j^{-|s_j|}{\sigma_j^{-n_j}}
\]
for non-zero integers $s_1,\ldots,s_k$ and $\sigma_j:=\mbox{signum}(s_j)$,
and
\[
{\tt zp}(p,s_1,\ldots,s_k)
  :=\sum_{n_1>\cdots>n_k}
    p^{-n_1}\prod_{j=1}^k n_j^{-s_j}
\]
for real $p\ge 1$ and positive integers $s_1,\ldots,s_k$.  The
code for evaluating these sums was written in C, using routines
from GMP, the GNU Multiprecision
Library\footnote{{\tt http://www.swox.com/gmp/}}.  Our implementation
permits the precision of the evaluation to be set anywhere between
10 and 100 digits.  Progress is currently underway
to extend the scope of sums that can be evaluated.  The exact status
of the EZ Face is at any moment documented at its ``Definitions''
and ``Using EZ-Face'' pages.

In addition to the functions {\tt z} and {\tt zp}, the {\tt
lindep} function, based on the LLL
algorithm~\cite{LLL} for discovering integer relations~\cite{BL}
satisfied by a vector of real numbers, can be called. An integer
relation for a vector of real numbers $(x_1,\ldots,x_n)$ is a
non-zero integer vector $(c_1,\ldots,c_n)$ such that $\sum_{i=1}^n
c_ix_i=0$. The required syntax is ${\tt
lindep}([x_1,\ldots,x_n])$, where $x_1,\ldots,x_n$ is the vector
of values for which the relation is sought.  One must ensure that
the vector of real numbers is evaluated to sufficient precision to
avoid bogus relations and other numerical artifacts.  The
{\tt lindep} code was written by Michael Monagan
and Greg Fee, both of the CECM,
and is available on request.  Send e-mail to either
{\tt monagan@cecm.sfu.ca} or {\tt gjfee@cecm.sfu.ca}.

%{\sc EZ Face Examples.}

Below, we give some examples showing how EZ Face may be used. The
left-aligned lines represent the input to EZ Face, while the
centered lines represent the output of EZ Face. All computations
are done with the precision of 50 digits.

\begin{example}\label{ex:z6}
\begin{verbatim}

Pi^6/z(6)

         945.00000000000000000000000000000000000000000000000
\end{verbatim}
\end{example}

\begin{example}\label{ex:z413-red}
\begin{verbatim}

lindep([z(4,1,3), z(5,3), z(8), z(5)*z(3), z(3)^2*z(2)])

                     36., 36., -71., 90., -18.
\end{verbatim}
\end{example}

\begin{example}\label{ex:z3-zetaLilog}
\begin{verbatim}

lindep([ z(3), Pi^2*log(2), zp(2,2,1), zp(2,3) ])

                       12., -1., -12., -12.
\end{verbatim}
\end{example}
Example~\ref{ex:z6} is a simple instance of Euler's formula for
$\z(2n)$. Example~\ref{ex:z413-red} is the discovery of
equation~(\ref{z413-red}).  Example~\ref{ex:z3-zetaLilog} confirms
formula (\ref{z3-zetaLilog}).

\section{Evaluations for Unit Euler Sums}
\label{S8}
As usual, the H\"older conjugates $p$ and $q$ denote real numbers
satisfying $1/p+1/q=1$, and $p>1$ or $p\le-1$ for convergence.
Our first result is an easy consequence of the binomial theorem.
\begin{theorem}
\label{T3}
The generating function equality
\[
   1+\sum_{n=1}^\infty x^n \mu(\{p\}^n) = q^x.
\]
holds.
\end{theorem}
\begin{proof} By definition~(\ref{unitEuler}) of $\mu$,
\begin{eqnarray*}
   1+\sum_{n=1}^\infty x^n \mu(\{p\}^n)
   &=& 1+x\sum_{m=1}^\infty \frac1{mp^m}\prod_{j=1}^{m-1}\(1+\frac{x}j\)\\
   &=& 1+\sum_{m=1}^\infty \(\frac{-1}p\)^m{-x\choose m}\\
   &=& \(1-1/p\)^{-x}\\
   &=& q^x.
\end{eqnarray*}
%\eop
\end{proof}

\begin{Cor}
\label{C1}
\[
   \mu(\{p\}^n) = (\log q)^n/n!,\qquad 0\le n\in\Z.
\]
\end{Cor}

\begin{remarks} Of course, when $n=0$, we need to invoke the
usual empty product convention to properly interpret $\mu(\{\})=1$.
Since the mapping $p\mapsto 1-p$ induces the mapping $q\mapsto 1/q$ under
the H\"older correspondence, duality~(\ref{simpledual}) takes the
particularly appealing form $\mu(\{p\}^n)=(-1)^n\mu(\{1-p\}^n)$
in this context.
In particular, $p=-1$ and
$\d$-duality~(\ref{deltadual}),~(\ref{deltaus}) gives
\[
   \d(\us^n) = (-1)^n\mu(\{-1\}^n) = (\log 2)^n/n!,
   \qquad 0\le n\in\Z,
\]
i.e.
\[
\begin{split}
   \sum_{\nu_1,\ldots,\nu_n=1}^\infty \;
   \prod_{j=1}^n\frac{1}{2^{\nu_j}(\nu_j+\cdots+\nu_n)}
   &=\sum_{\nu_1,\ldots,\nu_n=1}^\infty \;
   \prod_{j=1}^n\frac{(-1)^{\nu_j+1}}{\nu_j+\cdots+\nu_n}\\
   &=\frac{\(\log 2\)^n}{n!},\qquad 0\le n\in\Z,
\end{split}
\]
which can be viewed as an iterated sum extension of the well-known result
\[
   \sum_{\nu=1}^\infty\frac1{\nu 2^\nu}
   =\sum_{\nu=1}^\infty\frac{(-1)^{\nu+1}}\nu
   =\log 2,
\]
typically obtained by comparing the Maclaurin series for $\log(1+x)$ when
$x=-\tfrac12$ and $x=1$.
\end{remarks}

We now prove a~few results for unit Euler sums that were left as
open conjectures in~\cite{BBB}.  It will be convenient to employ
the following notation:
\begin{equation}
   A_r := {\Li}_r(\tfrac12) = \d(r) = \sum_{k=1}^\infty\frac1{2^k k^r},
   \quad P_r := \frac{(\log 2)^r}{r!},
   \quad Z_r := (-1)^r\z(r).
   \label{APZdefs}
\end{equation}
\begin{theorem}
\label{T4}
For all positive integers $m$,
\[
   \mu(\{-1\}^m,1) = (-1)^{m+1}\sum_{k=0}^m A_{k+1}P_{m-k} - Z_{m+1}.
\]
\end{theorem}

\begin{proof}
 From the case~(\ref{zetaLilog}) of the H\"older convolution, we have
\[
   \d(2,\us^{m-1}) = \z(m+1) - \sum_{r=1}^{m+1}\d(r)\d(\us^{m+1-r}).
\]
Now multiply both sides by $(-1)^m$ and apply the case~(\ref{delta2us})
of $\d$-duality.
%\eop
\end{proof}

\begin{remarks} Theorem~\ref{T4} appeared as the conjectured formula~(67)
in~\cite{BBB}, and is valid for all nonnegative integers $m$ if
the divergent $m=0$ case is interpreted appropriately.  The equivalent
generating function identity is
\[
\begin{split}
   \sum_{n=1}^\infty x^n\mu(\{-1\}^n,1)
   &=\int_0^{1/2} \frac{(1-t)^x-1}{t}dt\\
   &=\log 2+\sum_{n=1}^\infty\(\frac1{x+n}-\frac1n\)
   -\sum_{n=1}^\infty\frac{2^{-(x+n)}}{x+n},
\end{split}
\]
correcting the misprinted sign in formula~(21) of~\cite{BBB}.
\end{remarks}

The asymmetry which marrs Theorem~\ref{T4} is recovered in the
generalization~(\ref{BBB68}), restated and proved below.
\begin{theorem}
\label{T5}
For all positive integers $m$ and all nonnegative integers $n$,
we have
\begin{eqnarray}
   \mu(\{-1\}^m,1,\{-1\}^n) &=&
      (-1)^{m+1}\sum_{k=0}^m{n+k\choose n}A_{k+n+1}P_{m-k}\nonumber\\
     &+& (-1)^{n+1}\sum_{k=0}^n{m+k\choose m}Z_{k+m+1}P_{n-k},
\end{eqnarray}
where $A_r$, $P_r$ and $Z_r$ are as in~(\ref{APZdefs}).
\end{theorem}

\begin{proof}
Let $m$ be a~positive integer, and let $n$ be a~nonnegative integer.
We have
\begin{eqnarray*}
   \mu(\{-1\}^m,1,\{-1\}^n)
   &=& (-1)^{m+n+1}\iu\w_{-1}^m\w_1\int_0^y\w_{-1}^n\\
   &=& (-1)^{m+n+1}\iu\w_{-1}^m\w_1\int_1^{1-y}\w_2^n\\
   &=& (-1)^{m+n+1}\iu\w_{-1}^m\w_1\int_{1/2}^{(1-y)/2}\w_1^n\\
   &=& (-1)^{m+n+1}\iu\w_{-1}^m\w_1\(\log(1+y)\)^n/n!.
\end{eqnarray*}
By duality,
\begin{eqnarray*}
   m!n!\mu(\{-1\}^m,1,\{-1\}^n)
   &=& m!\iu\(-\log(2-y)\)^n\w_0\w_2^m\\
   &=& m!\iu\(-\log(2-y)\)^n\w_0\int_0^{y/2}\w_1^m\\
   &=&\iu\(-\log(2-y)\)^n\(\log(1-y/2)\)^m\,dy/y.\\
\end{eqnarray*}
Letting $t=1-y/2$ and forming the generating function, it follows that
\[
\begin{split}
   &\quad\sum_{m=1}^\infty\sum_{n=0}^\infty x^m
   y^n\mu(\{-1\}^m,1,\{-1\}^n)\\
   &= \sum_{m=1}^\infty\sum_{n=0}^\infty\frac{x^m}{m!}\frac{y^n}{n!}
      \int_{1/2}^1\(-\log(2t)\)^n\(\log t\)^m\frac{dt}{1-t}\\
   &= \int_{1/2}^1 \frac{(2t)^{-y}\(t^x-1\)}{1-t}\,dt.
\end{split}
\]
Expanding $1/(1-t)$ in powers of $t$ and integrating term by term
yields
\begin{eqnarray}
   &&\sum_{m=1}^\infty\sum_{n=0}^\infty
      x^m y^n\mu(\{-1\}^m,1,\{-1\}^n)\nonumber\\
   &&\quad= 2^{-y}\sum_{k=1}^\infty\(\frac1{k+x-y}-\frac1{k-y}\)
   -\sum_{k=1}^\infty\frac{2^{-(k+x)}}{k+x-y}
   +\sum_{k=1}^\infty\frac{2^{-k}}{k-y}.
   \label{muGF1}
\end{eqnarray}
%It is now a~routine matter to extract the coefficient of $x^m y^n$.
%For the record, here are the details:
Since $m\ge1$, we may ignore the terms in~(\ref{muGF1}) which are
independent of $x$.  Thus formally, but with the divergences
coming only from the terms independent of $x$ and hence harmless,
\begin{eqnarray*}
   &&   -2^{-x}\sum_{k=1}^\infty
   \frac{2^{-k}}{k+x-y}+2^{-y}\sum_{k=1}^\infty\frac{1}{k+x-y}\\
   && \qquad = -\sum_{r=0}^\infty (-x)^r P_r\sum_{h=1}^\infty
   (y-x)^{h-1}A_h-\sum_{r=0}^\infty (-y)^rP_r\sum_{h=1}^\infty
   (x-y)^{h-1}Z_h,
\end{eqnarray*}
where we have used the abbreviations in~(\ref{APZdefs}).  It is
now a~routine matter to extract the coefficient of $x^m y^n$ to
complete the proof.
%\eop
\end{proof}

\begin{remark}
Theorem~\ref{T5} is an extension of conjectured formula~(68) of~\cite{BBB},
and is valid for all nonnegative integers $m$ and $n$ if the divergent
$m=0$ case is interpreted appropriately.
\end{remark}

\section{Other Integral Transformations}
\label{sec-int-transf}

In Section~\ref{S6}, we proved the duality principle for multiple
polylogarithms by using the integral transformation $y\mapsto
1-x$.  Similarly, in this section we prove additional results for
multiple polylogarithms by using suitable transformations of
variables in their integral representations.

\begin{theorem}
\label{cyclotomicThm}
Let $n$ be a positive integer.  Let $b_1,\ldots,b_k$
be arbitrary complex numbers, and let $s_1,\ldots, s_k$ be positive
integers.  Then
\[
   \lambda{s_1,s_2,\ldots,s_k\choose b_1^n,b_2^n,\ldots,b_k^n}
= n^{s-k} \sum \lambda{s_1,\ldots,s_k\choose\eps_1 b_1,\ldots,\eps_k b_k},
\]
where the sum is over all $n^k$ cyclotomic sequences
\[
   \eps_1,\ldots,\eps_k \in \left\{1,e^{2\pi i/n}, e^{4\pi i/n},\ldots,
   e^{2\pi i(n-1)/n}\right\},
\]
and, as usual, $s:=s_1+s_2+\cdots + s_k$.
\end{theorem}

\begin{proof}  Write the left-hand side as an iterated integral
as in (\ref{iterint}):
\[
   L :=  \lambda{s_1,s_2,\ldots,s_k\choose b_1^n,b_2^n,\ldots,b_k^n}
      = (-1)^k \int_0^1 \prod_{j=1}^k \w_0^{s_j-1}\w(b_j^n).
\]
Now let $y=x^n$ at each level of integration.  This sends $\w_0$
to $n\w_0$ and, by partial fractions,
\[
   \w(b^n) \mapsto \sum_{r=0}^{n-1} \w\left(b e^{2\pi i r/n}\right).
\]
The change of variable gives
\[
   L = (-1)^k \int_0^1 \prod_{j=1}^k (n\w_0)^{s_j-1}
       \sum_{r=0}^{n-1} \w\left( b_j e^{2\pi i r/n} \right).
\]
Now carefully expand the noncommutative product and reinterpret
each resulting iterated integral as a $\l$-function
to complete the proof.
%\eop
\end{proof}

%{\sc Example.}
\begin{example}
When $n=2$ and $k=1$, Theorem~\ref{cyclotomicThm} asserts that
\[
   \z(s)=2^{s-1}\sum_{n=1}^\infty\frac{1+(-1)^n}{n^s}.
\]
Thus, Theorem~\ref{cyclotomicThm} can be viewed as a cyclotomic
extension of the well-known ``sum over signs'' formula for the
alternating zeta function:
\[
   \sum_{n=1}^\infty \frac{(-1)^{n+1}}{n^s}
   = (1-2^{1-s})\z(s),\quad \Re(s)>0.
\]
\end{example}

Next we prove two broad generalizations of
formulae (24), (26) and (28) of \cite{BBB}. By a pair of $\Cat{}{}$ operators
we mean nested concatenation
(similarly as two $\sum$ signs mean nested summation).

\begin{theorem}
\label{thm_la-1_to_mu}
Let $s_1,s_2,\ldots,s_k$ be nonnegative integers.
Then
%     lambda(1+s_k, 1+s_{k-1},..., 1+s_1)
%      ( -1,    -1,    ...,   -1   )
%
%     ---       k       s_j               ----
%   = >    mu( Cat {-1} Cat {eps_{i,j}} ) |  | eps_{i,j},
%     ---      j=1      i=1               |  |
%                                       1<=j<=k
%                      1<=i<=s_j
\[
     \l\left(\begin{array}{cccc} 1+s_k,& 1+s_{k-1},& \ldots,& 1+s_1 \\
                                 -1,&    -1,&    \ldots,&   -1
             \end{array} \right)
     =\sum\mu\bigg(\Cat{j=1}{k}\{-1\}\Cat{i=1}{s_j}\{\eps_{i,j}\}\bigg)
      \prod_{j=1}^k \prod_{i=1}^{s_j} \eps_{i,j}
\]
where the sum is over all $2^{s_1+s_2+\cdots+s_k}$ sequences of signs
$(\eps_{i,j})$, with each $\eps_{i,j}\in\{1,-1\}$ for all $1\le i\le s_j$,
$1\le j\le k$, and $\Cat{}{}$ denotes string concatenation.
\end{theorem}

\begin{proof}
Let
\[
   L :=    \l\left(\begin{array}{cccc}1+s_k,& 1+s_{k-1},& \ldots,& 1+s_1 \\
                                      -1,&    -1,&    \ldots,&   -1
                   \end{array} \right)
      = (-1)^k \int_0^1 \prod_{j=k}^1 \w_0^{s_j} \w_{-1}.
\]
Now let us use duality, and then we let $y = 2t/(1+t)$ at each level
of integration.  We get
\[
   L = (-1)^k \int_0^1 \prod_{j=1}^k \w_{-1} (\w_{-1}-\w_1)^{s_j}.
\]
Now let us carefully expand the noncommutative product.
We get
\[
   L = (-1)^k \sum (-1)^{\# \eps_{i,j}=1}
       \int_0^1 \prod_{j=1}^k \w_{-1} \prod_{i=1}^{s_j} \w(\eps_{i,j}),
\]
where the sum is over all sign choices $\eps_{i,j} \in \{1,-1\}$,
$1\le i\le s_j$, $1\le j\le k$,
and where by $\# \eps_{i,j}=a$ we mean the cardinality
of the set $\{(i,j)\;|\;\eps_{i,j}=a\}$.

Let us now interpret the iterated integrals as $\l$-functions.  In this
case, they are all unit Euler $\mu$-sums, as we defined in (\ref{unitEuler}).
Thus,
\[
   L = (-1)^k \sum (-1)^{\# \eps_{i,j}=1} (-1)^{k+s}
       \mu\bigg( \Cat{j=1}{k} \{-1\}
       \Cat{i=1}{s_j} \{\eps_{i,j}\} \bigg),
\]
where, as usual, $s:=s_1+s_2+\cdots+s_k$.  Now if $r$ of the $\eps_{i,j}$
equal $+1$, then $s-r$ of them equal $-1$.  Hence,
\[
   L = \sum (-1)^{\# \eps_{i,j}=-1}
    \mu\bigg( \Cat{j=1}{k} \{-1\} \Cat{i=1}{s_j} \{\eps_{i,j}\} \bigg).
\]
Finally, $(-1)^{\# \eps_{i,j}=-1}$ is the same as the product over
all the signs $\eps_{i,j}$, and this latter observation completes
the proof of Theorem~\ref{thm_la-1_to_mu}.
%\eop
\end{proof}

Theorem~\ref{thm_la-1_to_mu}
generalizes several identities conjectured in \cite{BBB}.
For example, we get the conjecture (28) of \cite{BBB}
if we put $s_{n+1}=m$, $s_n=s_{n-1}=...=s_1=0$ in Theorem~\ref{thm_la-1_to_mu}.
Furthermore, (24) of \cite{BBB} is the case $s_{m+n+1}=s_{m+n}=...=s_{n+2}=0$,
$s_{n+1}=1$, $s_n=s_{n-1}=...=s_1=0$, and (26) of \cite{BBB} is a special
case of Theorem~\ref{thm_la-1_to_mu} as well.
Thus {\em every} multiple polylogarithm with all alternations (or,
equivalently, every Euler sum with first position alternating and
all the others non-alternating) is a signed sum over unit Euler
sums. The representation of the sign coefficients used in
Theorem~\ref{thm_la-1_to_mu} is much simpler than the cumbersome
form of (28) in \cite{BBB}.

Below we present
a dual to Theorem~\ref{thm_la-1_to_mu}, which gives
{\em any} unit Euler $\mu$-value in terms of $\l$-values
with all alternations (equivalently,
Euler sums with only first position alternating):

\begin{theorem} Let $s_1,s_2,...,s_k$ be nonnegative integers.
Then
%    k-1
%      mu( \Cat {-1}{1}^{s_{k-j}} )
%    j=0
%
%
%      ---           k   q_j
%    = >    lambda( \Cat  \Cat t_{i,j} )
%      ---        ( j=1  i=1         )
%       ( -1, -1, ...., -1 )
\[
     \mu\bigg( \Cat{j=0}{k-1} \{-1\}\{1\}^{s_{k-j}} \bigg)
   = \sum \l\bigg(\Cat{j=1}{k}\Cat{i=1}{q_j}\{t_{i,j}-\}\bigg)
\]
where the sum is over all $2^{s_1+s_2+\cdots+s_k}$ positive integer
compositions
\[
   t_{1,j} + t_{2,j} + ... + t_{q_j, j} = s_j + 1,\qquad
   1\le q_j\le s_j+1,\ \ 1\le j\le k.
\]
\end{theorem}

\begin{proof}
Let
\[
   M :=  \mu\bigg( \Cat{j=0}{k-1} \{-1\}\{1\}^{s_{k-j}} \bigg)
   =  (-1)^k \delta\bigg( \Cat{j=1}{k} \{1+s_j\} \bigg)
   =  \int_0^1 \prod_{j=1}^k \w_0^{s_j} \w_2.
\]
Again, let us make the change of variable $y=2t/(1+t)$ at each level.  Then
\[
   M = \int_0^1 \prod_{j=1}^k (\w_0-\w_{-1})^{s_j} (-\w_{-1}).
\]
Again, let us carefully expand the noncommutative product.  We get
\[
   M = \sum (-1)^{\# \eps_{i,j}=-1} \int_0^1 \prod_{j=1}^k \left[
    \prod_{i=1}^{s_j} \w(\eps_{i,j}) \right] (-\w_{-1}),
\]
where this time, the sum is over all $\eps_{i,j}\in\{0,-1\}$ with
$1\le i\le s_j,\ 1\le j\le k$.
Note that each $\w_{-1}$ in the integrand contributes $-1$ to the sign
and $+1$ to the depth.  Since
\[
   (-1)^{\mbox{depth}} \int_0^1 \mbox{weight-length string}
   = \lambda(\mbox{depth-length string}),
\]
it follows that $M$ is a sum of $\l$-values with all $+1$ coefficients.
That is,
\[
   M = \sum \lambda{ \vec t_1, \ldots ,\vec t_k \choose -1,\ldots, -1},
\]
where the sum is over all vectors
\[
   \vec t_j = (t_{1,j}, \ldots, t_{q_j,j}),\qquad  1\le q_j\le 1+s_j,
\]
and such that
\[
   \sum_{i=1}^{q_j} t_{i,j} = 1+s_j,\qquad 1\le j\le k.
\]
In other words, the sum is over all $2^s$ independent
positive integer compositions (in the technical sense of
combinatorics) of the numbers $1+s_j$, $1\le j\le k$.
%\eop
\end{proof}

\section{Functional Equations}

One fruitful strategy for proving identities involving special
values of polylogarithms is to prove more general (functional, differential)
identities and instantiate them at appropriate argument values.
In the last two sections of this paper we present examples of such proofs.

\begin{lemma}
\label{lem-eq-A12}
Let $0\le x\le 1$ and let
\[
   J(x):=\int_0^x\frac{\left(\log(1-t)\right)^2}{2t}\,dt
\]
Then
\begin{eqnarray}
\label{lem-A12}
   J(-x) = -J(x)+\tfrac14J(x^2)+J\!\left(\frac{2x}{x+1}\right)
         -\tfrac18J\!\left(\frac{4x}{(x+1)^2}\right).
\end{eqnarray}
\end{lemma}

\begin{proof}
If $L(x)$ and $R(x)$
denote the left-hand and the right-hand sides of (\ref{lem-A12}),
respectively, then by elementary manipulations (under the assumption $0<x<1$)
we can show that $dL/dx=dR/dx$.
The easy observation $L(0)=R(0)=0$
then completes the proof.
%\eop
\end{proof}

\begin{remarks}
The identity (\ref{lem-A12}) can be discovered and proved using a
computer.  Once the ``ingredients'' (the $J$-terms) of the
identity are chosen, the constant coefficients at them can be
determined by evaluating the $J$-terms at a sufficiently arbitrary
value of $x\in ]0,1[$ and using an integer relation algorithm
\cite{BL}.  Once the identity is discovered, the main part of the
proof (namely showing that $dL/dx=dR/dx$) can be accomplished in a
computer algebra system (e.g., using the {\tt simplify()} command
of Maple).
\end{remarks}

\begin{theorem}
\label{thm21-21}
We have
\begin{equation}
\label{th21-21}
\l(2-,1-) = \z(2,1)/8.
\end{equation}
\end{theorem}

\begin{proof}
Using notation of Lemma~\ref{lem-eq-A12}
let us observe that
\[
   J(x)=\sum_{n_1>n_2>0} \frac{x^{n_1}}{n_1^2 n_2}.
\]
Plugging in $x=1$ and applying (\ref{lem-A12}) now completes the proof.
%\eop
\end{proof}

\begin{remarks}
Theorem~\ref{thm21-21} is the $n=1$ case of the conjectured
identity (23) of \cite{BBB}, namely
\begin{equation}
\label{eq21n-21n}
   \l(\underbrace{2-,1-,2,1,\ldots}_{2n})
   \eu 8^{-n}\z(\{2,1\}^n),
\end{equation}
for which we have
overwhelming numerical evidence.
This evidence also
suggests that (\ref{eq21n-21n}) with $n>1$ seems to be the only case
when two Euler sums that do not evaluate (in the sense
of the definition in Section~\ref{S3})
have a~rational quotient, different from 1. (See also Section~\ref{S62}.)
\end{remarks}

\section{Differential Equations and Hypergeometric Series}
\label{S10}
Here, it is better to work with
\[
   L(s_1,\ldots,s_k;x) := \l_{1/x}(s_1,\ldots,s_k),
\]
since then we have
\[
   \frac{d}{dx} L(s_k,\ldots,s_1;x) = \frac1{x}L(-1+s_k,\ldots,s_1;x)
\]
if $s_k\ge 2$; while for $s_k=1$,
\[
   \frac{d}{dx} L(s_k,\ldots,s_1;x) = \frac{1}{1-x}L(s_{k-1},\ldots,s_1;x).
\]
With the initial conditions
\[
   L(s_k,\ldots,s_1;0) = 0,\quad k\ge 1,
   \quad{\mbox{and}}\quad L(\{\};x) :=1,
\]
the differential equations above determine the $L$-functions uniquely.

\subsection{Periodic Polylogarithms}
\label{S101}
If $\vec s := (s_1,s_2,\ldots,s_k)$ and $s:=\sum_{j}s_j$,
then every {\em periodic polylogarithm} $L(\{\vec s\}^r)$ has an ordinary
generating function
\[
   L_{\vec s}(x,t) := \sum_{r=0}^\infty L(\{\vec s\}^r;x) t^{rs}
\]
which satisfies an algebraic ordinary differential equation in $x$.
In the simplest case, $k=1$, $\vec s$ reduces to the scalar
$s$, and the differential equation for the ordinary generating
function is $D_s - t^s =0$, where
\[
   D_s := \((1-x)\frac{d}{dx}\)^1\(x\frac{d}{dx}\)^{s-1}.
\]
The series solution is a~generalized hypergeometric function
\begin{eqnarray*}
   L_{s}(x,t)
   &=& 1+\sum_{r=1}^\infty x^r\frac{t^s}{r^s}\prod_{j=1}^{r-1}
   \(1+\frac{t^s}{j^s}\)\\
   &=& {}_sF_{s-1}\(\begin{array}{c}-\w t,-\w^3t,\ldots,-\w^{2s-1}t\\
   1,1,\ldots,1\end{array}
   \left|\phantom{\displaystyle\int}\!\!\!\!x\right.\),
\end{eqnarray*}
where $\w=e^{\pi i/s}$, a~primitive $s$th root of $-1$.

\subsection{Proof of Zagier's Conjecture}
\label{S102}
Let ${}_2F_1(a,b;c; x)$ denote the Gaussian hypergeometric
function.  Then:

\begin{theorem}
\label{T7}
\begin{multline}
   \sum_{n=0}^\infty L(\{3,1\}^n; x)t^{4n}\\
   ={}_2F_1\(\tfrac12 t(1+i), -\tfrac12 t(1+i); 1; x\)
    {}_2F_1\(\tfrac12 t(1-i), -\tfrac12 t(1-i); 1; x\).
    \label{ZAGsplit}
\end{multline}
\end{theorem}

\begin{proof}
Both sides of the putative identity start
\[
   1+\frac{t^4}8 x^2 + \frac{t^4}{18}x^3 + \frac{t^8+44t^4}{1536}x^4+\cdots
\]
and are annihilated by the differential operator
\[
   D_{31}:=\((1-x)\frac{d}{dx}\)^2\(x\frac{d}{dx}\)^2-t^4.
\]
Once discovered, this can be checked in Mathematica or Maple.
%\eop
\end{proof}

\begin{Cor} {\rm{(Zagier's Conjecture)\cite{Zag}}}
\label{C2}
For all nonnegative integers $n$,
\[
   \z(\{3,1\}^n) = \frac{2\pi^{4n}}{(4n+2)!}.
\]
\end{Cor}

\begin{proof} Gauss's ${}_2F_1$ summation theorem gives
\[
   {}_2F_1(a,-a;1;1) = \frac{1}{\G(1-a)\G(1+a)} = \frac{\sin(\pi a)}{\pi a}.
\]
Hence, setting $x=1$ in the generating function~(\ref{ZAGsplit}), we have
\[
\begin{split}
   &\quad\;\sum_{n=0}^\infty \z(\{3,1\}^n)t^{4n}\\
   &={}_2F_1\(\tfrac12 t(1+i), -\tfrac12 t(1+i); 1; 1\)
   {}_2F_1\(\tfrac12 t(1-i), -\tfrac12 t(1-i); 1; 1\)\\
   &=\frac{2\sin(\tfrac12 (1+i)\pi t)
           \sin(\tfrac12(1-i)\pi t)}{\pi^2 t^2}\\
   &=\frac{\cosh(\pi t)-\cos(\pi t)}{\pi^2 t^2}\\
   &=\sum_{n=0}^\infty \frac{2\pi^{4n} t^{4n}}{(4n+2)!}.
\end{split}
\]
%\eop
\end{proof}

\begin{remark}
The proof is Zagier's modification of Broadhurst's, based on the
extensive empirical work begun in~\cite{BBB}.
\end{remark}

\subsection{Generalizations of Zagier's Conjecture}
In \cite{BBBLc} we give an alternative (combinatorial) proof of
Zagier's conjecture, based on combinatorial manipulations of the
iterated integral representations of MZVs (see Sections~\ref{S42}
and~\ref{S53}). Using the same technique, we prove in \cite{BBBLc}
the
 ``Zagier dressed with~2'' identity:
\begin{equation}
\label{zagier-dressed}
\sum_{\vec s} \zeta(\vec s) = {{\pi^{4n+2}}\over{(4n+3)!}}
\end{equation}
where $\vec s$ runs over all $2n+1$ possible insertions of the number 2
in the string $\{3,1\}^n$.
Still,
(\ref{zagier-dressed}) is just the beginning of a large family
of conjectured identities that we discuss in \cite{BBBLc}.

\section{Open Conjectures}

The reader has probably noticed that many formulae proved in
this paper were conjectured in~\cite{BBB}.  For the sake
of completeness, we now list formulae from~\cite{BBB}
that are still
open:  (18), (23), (25), (27), (29), (44), and (70)--(74).
It is possible that some of these conjectures can be proved
using techniques of the present paper.

\section*{Acknowledgements}
%\acknowledgements{
Thanks are due to David Borwein, Douglas Bowman
and Keith Johnson for their helpful comments.  We are especially
grateful to the referee for a thorough and detailed report which
led to several improvements.
%}

\end{document}